\documentclass[11pt,a4paper]{article}
\usepackage[top=24mm,bottom=25mm,left=24mm,right=24mm,headheight=16pt]{geometry}
\usepackage{amsmath,amsthm,fontspec,unicode-math}
\usepackage{booktabs,longtable,array,enumitem,fancyhdr,lastpage,xcolor,xurl,needspace}

\definecolor{inkblue}{RGB}{31,57,77}
\usepackage[unicode,colorlinks=true,linkcolor=inkblue,citecolor=inkblue,urlcolor=inkblue,bookmarksnumbered=true]{hyperref}
\hypersetup{pdftitle={A Goldbach-Type Problem for the Unitary Möbius Function},pdfsubject={Pointwise lower bounds, explicit existence, and average asymptotics},pdfauthor={Tian Caolin and Yang Chenxi},colorlinks=true,
    urlcolor=blue,
    linkcolor=blue}
\usepackage{titlesec}\titleformat{\section}{\Large\bfseries\color{inkblue}}{\thesection}{1em}{}\titleformat{\subsection}{\large\bfseries}{\thesubsection}{1em}{}
\setlist[enumerate]{leftmargin=2em,itemsep=0.18em,topsep=0.3em}
\numberwithin{equation}{section}
\newtheoremstyle{paper}{6pt}{6pt}{\normalfont}{}{\bfseries}{.}{0.5em}{}
\theoremstyle{paper}
\newtheorem{theorem}{Theorem}[section]
\newtheorem{lemma}[theorem]{Lemma}
\newtheorem{proposition}[theorem]{Proposition}
\newtheorem{conjecture}[theorem]{Conjecture}
\newtheorem{corollary}[theorem]{Corollary}
\newtheorem{remark}[theorem]{Remark}
\newtheorem{definition}[theorem]{Definition}

\newcommand{\one}{\mathbf{1}}
\newcommand{\CB}{\mathfrak{C}_B}

\newcommand{\R}{\mathbb R}
\newcommand{\sgnset}{\{-1,+1\}}
\DeclareMathOperator{\rad}{rad}
\newcommand{\e}{\mathrm e}
\fancypagestyle{plain}{\fancyhf{}\fancyfoot[C]{\small \thepage\,/\,\pageref*{LastPage}}}
\title{\vspace{-1em}\bfseries A Goldbach-Type Problem for the Unitary Möbius Function}
\author{Tian Caolin\textsuperscript{1,2}\thanks{Email: \href{mailto:tiancaolin@lsec.cc.ac.cn}{\nolinkurl{tiancaolin@lsec.cc.ac.cn}}.}\quad Yang Chenxi\textsuperscript{1,2}\thanks{Email: \href{mailto:yangchenxi@lsec.cc.ac.cn}{\nolinkurl{yangchenxi@lsec.cc.ac.cn}}.}\\[0.6em]\parbox{0.92\textwidth}{\centering\small\textsuperscript{1}Institute of Computational Mathematics and Scientific/Engineering Computing,\\Academy of Mathematics and Systems Science, Chinese Academy of Sciences\\[0.35em]\textsuperscript{2}School of Mathematical Sciences, University of Chinese Academy of Sciences}}
\date{November 2, 2026}
\begin{document}
\maketitle
\begin{abstract}
Let $g(n)=(-1)^{\omega(n)}$, where $\omega(n)$ denotes the number of distinct prime factors of $n$. We prove that, for every sufficiently large integer $N$, each ordered additive sign pattern $(g(a),g(N-a))=(s,t)$ occurs at least $N/(128(\log N)^3)$ times. We also prove that, for every $N\ge18$ and every $s,t\in\sgnset$, there is at least one integer $1\le a<N$ such that $(g(a),g(N-a))=(s,t)$. In addition, we establish a one-quarter asymptotic formula for almost all target integers. Finally, we discuss Goldbach-type problems for the ordinary Möbius function and formulate related questions.
\end{abstract}
\clearpage
\begingroup
\small
\linespread{1}\selectfont
\setlength{\parskip}{0pt}
\tableofcontents
\endgroup
\clearpage

\section{Introduction}\label{sec:introduction}

The Goldbach problem originated in the correspondence between Goldbach and Euler in 1742. Its strong form is usually stated as follows.
\begin{conjecture}[Strong Goldbach conjecture]
Every even integer greater than $2$ is the sum of two primes; equivalently, every even integer greater than $4$ is the sum of two odd primes.
\end{conjecture}
This conjecture remains unresolved. The weak Goldbach problem concerns representations of odd integers as sums of three primes. In 1937, Vinogradov proved the assertion for all sufficiently large odd integers~\cite{Vinogradov1937}.
\begin{theorem}[Vinogradov\textquotesingle s three-primes theorem]
There is a constant $N_0$ such that every odd integer $N$ greater than $N_0$ can be written as
\[
N=p_1+p_2+p_3,
\]
where $p_1,p_2,p_3$ are odd primes.
\end{theorem}
The full weak Goldbach conjecture, asserting that every odd integer greater than $5$ is the sum of three primes, was proved by Helfgott~\cite{Helfgott2013,HelfgottICM2014}.

These problems can also be formulated in terms of sumsets. Given subsets $A,B,C$ of the natural numbers, write
\[
A+B=\{a+b:a\in A,\ b\in B\}.
\]
The question is under what conditions $C\subseteq A+B$, or whether there exists $N_0$ such that
$C\cap(N_0,\infty)\subseteq A+B$. We refer to such additive representation problems collectively as Goldbach-type problems.
For example, if $\mathcal P$ is the set of primes and $C=\{4,6,8,\ldots\}$, then $C\subseteq\mathcal P+\mathcal P$ is the strong Goldbach conjecture.
Chen Jingrun\textquotesingle s $1+2$ theorem is a major advance in this direction~\cite{Chen1973}. In sumset notation, it can be stated as follows.
\begin{theorem}[Chen\textquotesingle s theorem]
Let
\[
B_2=\mathcal P\cup\{p_1p_2:p_1,p_2\in\mathcal P\}.
\]
Then there is a constant $N_0$ such that
\[
C\cap(N_0,\infty)\subseteq\mathcal P+B_2.
\]
In other words, every sufficiently large even integer is the sum of a prime and an integer with at most two prime factors, counted with multiplicity.
\end{theorem}

In their 2026 preprint~\cite{LiLiu2026}, Jiamin Li and Jianya Liu obtained a $1+1.9$ result. Their Theorem 1.1 can be stated as follows.
\begin{theorem}[The $1+1.9$ result in the Li--Liu preprint]
Let
\[
B_{1.9}=\mathcal P\cup
\{p_1p_2:p_1,p_2\in\mathcal P,\ p_1\le p_2^{0.9}\}.
\]
Then there is a constant $N_0$ such that
\[
C\cap(N_0,\infty)\subseteq\mathcal P+B_{1.9}.
\]
\end{theorem}
These existence results for sufficiently large integers do not automatically yield a small starting point suitable for exhaustive verification; an explicit threshold requires a separate estimate.

Let $\mathbb N$ denote the positive integers, and let $\omega(n)$ count the distinct prime factors of $n\in\mathbb N$, with $\omega(1)=0$. We study additive representations using the multiplicative function
\[
g(n)=(-1)^{\omega(n)}.
\]
This is usually called the unitary M\"obius function; for classical work on arithmetic functions associated with unitary divisors, see~\cite{Cohen1960}. For squarefree $n$, we have $g(n)=\mu(n)$; at nonsquarefree integers, the two functions generally differ.
For $s\in\sgnset$, write
\[
A_s=\{a\in\mathbb N:g(a)=s\}.
\]
For an integer $N\ge2$ and $s,t\in\sgnset$, define the ordered representation count
\begin{equation}\label{intro:R}
R_{st}(N)=\#\{1\le a<N:g(a)=s,\ g(N-a)=t\}.
\end{equation}
Both summands are positive integers. Unless $a=N-a$, exchanging the two summands produces a different ordered representation. Also define
\begin{equation}\label{intro:SC}
S_g(x)=\sum_{n\le x}g(n),\qquad
C_g(N)=\sum_{a=1}^{N-1}g(a)g(N-a).
\end{equation}
Throughout, the upper summation limit in $S_g$ is $n\le x$; the corresponding single sum in the pattern-counting identity is $S_g(N-1)$.

Since $g(p)=-1$ for every prime $p$, the strong Goldbach conjecture implies the existence of the negative-negative pattern for every even $N>2$, that is, $R_{--}(N)>0$.
However, $A_-$ also contains many composite integers; for instance, $g(30)=-1$. We therefore consider all sufficiently large integers, rather than only even integers. Our first main result is as follows.
\begin{theorem}[Quantitative lower bounds for all sufficiently large integers]\label{thm:main-log}
There is an absolute constant $N_1$ such that every integer $N\ge N_1$ satisfies
\begin{align}
R_{--}(N),R_{++}(N)&\ge\frac{N}{128(\log N)^3},\label{intro:same}\\
R_{-+}(N)=R_{+-}(N)&\ge\frac{N}{64(\log N)^3}.\label{intro:mixed}
\end{align}
\end{theorem}
In particular, $N\in A_s+A_t$ for every $N\ge N_1$ and $s,t\in\sgnset$. The value of $N_1$ obtained here depends on constants in classical mean-value estimates and has not been evaluated numerically. Separately, finite integer verification yields the following existence result.
\begin{theorem}\label{thm:main-explicit}
For integers $N\ge2$, the complete exceptional sets for the four ordered sign patterns are
\begin{align*}
R_{--}(N)=0&\iff N\in\{2,3\},\\
R_{-+}(N)=R_{+-}(N)=0&\iff N\in\{2,7\},\\
R_{++}(N)=0&\iff N\in\{3,4,5,6,8,9,10,14,17\}.
\end{align*}
Thus all four patterns occur simultaneously if and only if
\[
N\notin\{2,3,4,5,6,7,8,9,10,14,17\}.
\]
In particular, $R_{st}(N)>0$ for every $N\ge18$ and all $s,t\in\sgnset$.
This common eventual threshold is optimal, since $R_{++}(17)=0$. The optimal eventual thresholds for the negative-negative, the two mixed, and the positive-positive patterns are $4,8,18$, respectively.
\end{theorem}
In sumset language, Theorem~\ref{thm:main-explicit} asserts that every $N\ge18$ belongs to all four sets $A_s+A_t$.
The proof first treats all $N\ge100000$ using covering chains, finite witness tables, and an argument for the remaining infinite tail. It then checks $2\le N<100000$ individually to determine all exceptions.
The first part handles the infinite range and the second handles the finite range; together they give the complete classification and optimal thresholds above.
The lower bounds in Theorem~\ref{thm:main-log} may fall short of the true order of magnitude. We therefore ask whether, as $N\to\infty$ through all integers,
\[
R_{st}(N)=\frac N4+o(N).
\]
We do not prove this pointwise asymptotic formula, but establish it for almost all integers.
\begin{theorem}[A one-quarter asymptotic formula for almost all integers]\label{thm:main-average}
For any fixed $A,B>0$, apart from at most $O_{A,B}(X/(\log X)^A)$ integers $2\le N\le X$, all four patterns simultaneously satisfy
\[
R_{st}(N)=\frac N4+O_B\!\left(\frac{N}{(\log N)^B}\right).
\]
\end{theorem}

Goldbach-type problems with sign constraints have been studied previously. In 2018, a version involving the Liouville function was posed on MathOverflow~\cite{LiouvilleMO}.
The Liouville function is $\lambda(n)=(-1)^{\Omega(n)}$, where $\Omega$ counts prime factors with multiplicity. It is completely multiplicative: the identity
$\lambda(mn)=\lambda(m)\lambda(n)$ holds for all positive integers $m,n$, unlike the corresponding identity for $g$.

Mangerel~\cite[Theorem 1.2]{Mangerel2024} proved that
\[
\left|\sum_{a=1}^{N-1}\lambda(a)\lambda(N-a)\right|<N-1
\qquad(N\ge11),
\]
and that equality for $N\ge2$ holds precisely at $\{2,3,5,10\}$. This guarantees the mixed patterns and at least one same-sign pattern. The fact that the correlation sum is not extremal does not by itself guarantee both same-sign patterns separately.
His subsequent preprint~\cite[Theorem 1.2 and Proposition 1.4]{MangerelGRH} proves, under the generalized Riemann hypothesis for Dirichlet $L$-functions, that every sufficiently large even integer admits a negative-negative representation. Under the same hypothesis, for sufficiently large \emph{prime} $N$, each of the four Liouville sign patterns occurs
\[
\gg N\exp\{-C(\log\log N)^6\}
\]
times.

In his 2026 preprint~\cite{Meng}, Meng gives a qualitative existence result: all four Liouville sign patterns occur for every $N\ge2$ with $N\notin\{2,3,4,5,6,9,10\}$. The author states that AI tools were used to assist the research. His rigidity method uses scaling laws for completely multiplicative functions and cannot be transferred directly to $g$; for the latter, scaling requires checking, factor by factor, whether a prime factor is already present or newly introduced.
For recent work on related additive convolutions and average estimates, see also~\cite{CGZ,Krishnamoorthy}.
On the other hand, the positive-sign sets of both functions have natural density $1/2$:
\[
\lim_{x\to\infty}\frac{\#\{n\le x:\lambda(n)=1\}}{x}=\frac12,
\qquad
\lim_{x\to\infty}\frac{\#\{n\le x:g(n)=1\}}{x}=\frac12.
\]
The negative-sign sets likewise have density $1/2$. These density statements follow from classical Möbius mean-value estimates and the estimates in Section~\ref{sec:basic}.\footnote{For $\lambda$, use the identity $\lambda(n)=\sum_{d^2\mid n}\mu(n/d^2)$ and, for fixed $D$, split the sum into $d\le D$ and $d>D$. The first part is $o(x)$ by $M_\mu(y)=o(y)$, while the absolute value of the second is at most $x\sum_{d>D}d^{-2}$. Hence $\sum_{n\le x}\lambda(n)=o(x)$. For $g$, apply Proposition~\ref{bas:mean}.}
An independent-random-sign heuristic gives the same first-order density background for the two functions. For classical work on additive properties of random integer sequences, see Erd\H{o}s--R\'enyi~\cite{ErdosRenyi1960}.
In the fixed-probability model, each positive integer belongs to a random set $A$ independently with probability $p\in(0,1]$. The natural density of $A$ equals $p$ with probability $1$. Moreover, with probability $1$, there exists a threshold $N_0=N_0(A)$, depending on the random set $A$, such that every integer $N>N_0$ belongs to $A+A$.
The latter assertion has a direct proof. For fixed $N\ge3$, consider the $\lfloor(N-1)/2\rfloor$ disjoint pairs $(a,N-a)$ with $1\le a<N/2$.
Both integers in each pair belong to $A$ with probability $p^2$, and these events are independent. Hence
\[
\mathbb P(N\notin A+A)\le(1-p^2)^{\lfloor(N-1)/2\rfloor}.
\]
The right-hand side is summable over $N$, so the Borel--Cantelli lemma proves the assertion. As $p$ increases, this upper bound for the probability of a missing sum decays faster.
An arbitrary deterministic set of positive density need not have this property: for example, $A=2\mathbb N$ has density $1/2$, and $A+A$ still has density $1/2$.
Thus equal densities alone do not imply that the additive problems for these particular arithmetic functions have the same difficulty; the probabilistic model provides a heuristic.

By contrast, the two nonzero-sign sets of the Möbius function each have density
\[
\lim_{x\to\infty}\frac{\#\{n\le x:\mu(n)=1\}}{x}
=\lim_{x\to\infty}\frac{\#\{n\le x:\mu(n)=-1\}}{x}
=\frac3{\pi^2}<\frac12.
\]
This follows from $\sum_{n\le x}\mu(n)=o(x)$,
$\sum_{n\le x}\mu^2(n)=6x/\pi^2+O(\sqrt x)$, and the identity
$\mathbf1_{\{\mu(n)=\pm1\}}=(\mu^2(n)\pm\mu(n))/2$.
The squarefree counting formula follows directly by summing $\mu^2(n)=\sum_{d^2\mid n}\mu(d)$ and using
$\sum_{d\ge1}\mu(d)/d^2=1/\zeta(2)=6/\pi^2$.
The Möbius function also takes the value zero, and its zero set has natural density $1-6/\pi^2>3/\pi^2$. In our argument, the five patterns involving a zero value already admit unconditional asymptotic estimates, whereas the pointwise asymptotics for the four nonzero patterns still depend on an unresolved correlation conjecture; see Section~\ref{sec:mobius}.

The prime set and the set of integers with at most two prime factors appearing in Chen\textquotesingle s theorem both have natural density zero, reflecting the sparser additive setting of the original Goldbach problem. These comparisons provide motivation rather than a rigorous ranking of difficulty. Intuitively, sets of positive density might seem more likely to satisfy a Goldbach property, but this need not be so.
One can also partition other arithmetic functions into level sets and ask whether they satisfy the corresponding Goldbach property, or even construct random sets $A,B$ and study the sumset $A+B$ or the product set $A\cdot B:=\{ab:a\in A,\ b\in B\}$. Section~\ref{sec:mobius} will further discuss the nine sign patterns for the ordinary Möbius function.

The remainder of the paper is organized as follows. In Section~\ref{sec:qplus}, we give four decompositions of an integer, which yield mixed-sign representations whenever an admissible parameter exists, and prove Theorems~\ref{thm:main-log} and~\ref{thm:main-average}. Section~\ref{sec:explicit} proves Theorem~\ref{thm:main-explicit}. The formalization of the proof of Theorem~\ref{thm:main-explicit} was completed by Yang Chenxi; the code is available in the \href{https://github.com/ilsunshine/Generalized_Goldbach_conjecture}{GitHub repository}.

Unless stated otherwise, $\log$ denotes the natural logarithm and $p,q$ denote primes. We write $\rad(n)$ for the product of the distinct prime factors of $n$, and $(a,b)$ and $[a,b]$ for the greatest common divisor and least common multiple, respectively.


\section{Basic identities and mean-value estimates}\label{sec:basic}
We first give multiplication and scaling identities, then establish the pattern-counting formula and the mean-value estimates needed below.
\subsection{Multiplication identities and pattern counts}

\begin{lemma}\label{bas:scaling}
For any positive integers $m,n$,
\begin{equation}\label{bas:gcd}
g(mn)=g(m)g(n)g((m,n)).
\end{equation}
In particular, $g$ is multiplicative. For any prime $p$, integer $k\ge1$, and $n\ge1$,
\begin{equation}\label{bas:prime-scale}
g(p^kn)=
\begin{cases}-g(n),&p\nmid n,\\g(n),&p\mid n.\end{cases}
\end{equation}
\end{lemma}
\begin{proof}
Inclusion--exclusion for the sets of prime factors gives
\[
\omega(mn)=\omega(m)+\omega(n)-\omega((m,n)).
\]
Since $(-1)^{-r}=(-1)^r$ for every integer $r$, raising $-1$ to the powers on both sides gives~\eqref{bas:gcd}. If $(m,n)=1$, then $g((m,n))=g(1)=1$.
If $p\nmid n$, multiplication by $p^k$ introduces exactly one distinct prime factor; if $p\mid n$, the prime-factor set is unchanged. This proves~\eqref{bas:prime-scale}, regardless of the parity of $k$.
\end{proof}

Equation~\eqref{bas:gcd} provides the correction factor in the noncoprime case. Subsequent scaling arguments use~\eqref{bas:prime-scale}, distinguishing whether a prime factor is already present.

\begin{lemma}\label{bas:four-counts}
For every $N\ge2$ and $s,t\in\sgnset$,
\begin{equation}\label{bas:counts}
4R_{st}(N)=N-1+(s+t)S_g(N-1)+stC_g(N).
\end{equation}
Consequently,
\begin{align}
4R_{--}(N)&=N-1-2S_g(N-1)+C_g(N),\\
4R_{++}(N)&=N-1+2S_g(N-1)+C_g(N),\\
4R_{-+}(N)&=4R_{+-}(N)=N-1-C_g(N).
\end{align}
\end{lemma}
\begin{proof}
Since $g(n)\in\sgnset$, we have $\mathbf1_{g(n)=s}=(1+sg(n))/2$. Multiply the two indicator functions and sum over $1\le a<N$. The substitution $b=N-a$ gives $\sum_{a=1}^{N-1}g(N-a)=S_g(N-1)$, proving the formula.
\end{proof}

Write
\[
P_\eta(N)=\#\{1\le a<N:g(a)g(N-a)=\eta\}\qquad(\eta\in\sgnset).
\]
Then
\begin{equation}\label{bas:product-counts}
P_+(N)+P_-(N)=N-1,\qquad
N-1\pm C_g(N)=2P_\pm(N).
\end{equation}
In particular, $|C_g(N)|<N-1$ guarantees only that both positive and negative correlation products occur; it does not specify which same-sign pattern occurs.
We also have
\begin{equation}\label{bas:min-count}
\min_{s,t}R_{st}(N)\ge
\frac{N-1-|C_g(N)|-2|S_g(N-1)|}{4}.
\end{equation}

Let $*$ denote Dirichlet convolution, and let $\kappa$ denote the convolution kernel defined below.
\begin{lemma}\label{bas:kernel}
Define the multiplicative function $\kappa$ by
\[
\kappa(1)=1,\qquad \kappa(p^a)=1-a\quad(a\ge1).
\]
Then
\begin{equation}\label{bas:convolution}
g=\mu*\kappa.
\end{equation}
For every fixed $\sigma>1/2$,
\begin{equation}\label{bas:absolute}
H_\sigma:=\sum_{d\ge1}\frac{|\kappa(d)|}{d^\sigma}
=\prod_p\left(1+\frac{p^{-2\sigma}}{(1-p^{-\sigma})^2}\right)<\infty.
\end{equation}
\end{lemma}
\begin{proof}
The Dirichlet convolution $\mu*\kappa$ is multiplicative and equals $1$ at $1$. At a prime $p$ it equals $\kappa(p)-\kappa(1)=-1$; at a prime power $p^a$ with $a\ge2$ it equals
\[
\kappa(p^a)-\kappa(p^{a-1})=(1-a)-(2-a)=-1.
\]
Here we have used $\mu(p^j)=0$ for $j\ge2$. These values agree with $g(p^a)$, so multiplicativity proves~\eqref{bas:convolution} at every positive integer.

Put $z=p^{-\sigma}\in(0,1)$. Since $\kappa(p)=0$,
\[
\sum_{a\ge0}|\kappa(p^a)|z^a
=1+\sum_{a\ge2}(a-1)z^a
=1+\frac{z^2}{(1-z)^2}.
\]
For fixed $\sigma$, the local factor differs from $1$ by $O_\sigma(p^{-2\sigma})$. Since $2\sigma>1$, these differences are summable. Expanding the Euler product over finite sets of primes and then applying monotone convergence to the nonnegative terms proves both the equality and convergence in~\eqref{bas:absolute}.
\end{proof}

\begin{proposition}[Dirichlet series]\label{bas:dirichlet}
For $\Re s>1$, the following identities hold with absolute convergence:
\[
 \sum_{n\ge1}\frac{g(n)}{n^s}
 =\prod_p\frac{1-2p^{-s}}{1-p^{-s}}
 =\frac{H(s)}{\zeta(s)},\qquad
 H(s)=\sum_{n\ge1}\frac{\kappa(n)}{n^s}
 =\prod_p\frac{1-2p^{-s}}{(1-p^{-s})^2}.
\]
The series and Euler product defining $H(s)$ converge absolutely and locally uniformly in $\Re s>1/2$,
and hence define an analytic function. We do not assert that $H(s)$ is everywhere nonzero.
\end{proposition}
\begin{proof}
For $\Re s>1$, the identity $|g(n)|=1$ permits an absolutely convergent Euler-product expansion. Put $z=p^{-s}$; then
$1+\sum_{a\ge1}g(p^a)z^a=(1-2z)/(1-z)$.
On the other hand, $\kappa(p^a)=1-a$ gives
$1+\sum_{a\ge1}(1-a)z^a=(1-2z)/(1-z)^2$.
Multiplication by $\zeta(s)^{-1}=\prod_p(1-p^{-s})$ yields the result.
For any compact set $K\subset\{\Re s>1/2\}$, choose
$1/2<\sigma_0\le\inf_{s\in K}\Re s$. Lemma~\ref{bas:kernel} gives
$H_{\sigma_0}<\infty$, providing a uniform absolutely summable majorant for the series. The difference between each local factor and $1$
is $-p^{-2s}/(1-p^{-s})^2$, and the sum of its absolute values also converges uniformly.
\end{proof}

\subsection{Mean values and densities of the sign sets}
We also use the classical estimate for the Möbius summatory function (see, for example, Green--Tao~\cite[p.~542, before (1.1)]{GT}): there is an absolute constant $c>0$ such that
\begin{equation}\label{bas:mobius-classical}
M_\mu(x):=\sum_{n\le x}\mu(n)
\ll x\exp(-c\sqrt{\log x})\qquad(x\ge3).
\end{equation}
This yields the following result.
\begin{proposition}\label{bas:mean}
There is an absolute constant $c_0>0$ such that
\begin{equation}\label{bas:mean-exp}
S_g(x)\ll x\exp(-c_0\sqrt{\log x})\qquad(x\ge3).
\end{equation}
Consequently, for every fixed $A>0$,
\begin{equation}\label{bas:mean-log}
S_g(x)\ll_A\frac{x}{(\log x)^A}\qquad(x\ge3).
\end{equation}
\end{proposition}
\begin{proof}
The convolution identity and an interchange of finite sums give
\[
S_g(x)=\sum_{d\le x}\kappa(d)M_\mu(x/d).
\]
First suppose $x\ge9$. For $d\le\sqrt x$, we have $x/d\ge3$ and $\log(x/d)\ge\frac12\log x$, so
\[
\sum_{d\le\sqrt x}|\kappa(d)M_\mu(x/d)|
\ll H_1x\exp\left(-\frac c{\sqrt2}\sqrt{\log x}\right).
\]
For $d>\sqrt x$, the bound $|M_\mu(y)|\le y$ for $y\ge1$ and~\eqref{bas:absolute} give
\[
\sum_{\sqrt x<d\le x}|\kappa(d)M_\mu(x/d)|
\le x\sum_{d>\sqrt x}\frac{|\kappa(d)|}{d}
\le H_{3/4}x^{7/8}.
\]
Fix $0<c_0\le c/\sqrt2$. Since $x^{7/8}=O(x\exp(-c_0\sqrt{\log x}))$, we obtain~\eqref{bas:mean-exp}. Adjusting the constant includes $3\le x<9$. Finally, $\exp(-c_0\sqrt{\log x})\ll_A(\log x)^{-A}$ gives~\eqref{bas:mean-log}.
\end{proof}

\begin{corollary}[Densities of sign sets and pointwise asymptotics]\label{bas:density}
For $\varepsilon\in\sgnset$,
\[
 \#\{n\le x:g(n)=\varepsilon\}
 =\frac{\lfloor x\rfloor+\varepsilon S_g(x)}2=\frac x2+o(x).
\]
Moreover, $C_g(N)=o(N)$ is equivalent to $R_{st}(N)=N/4+o(N)$ for all four patterns simultaneously.
\end{corollary}
\begin{proof}
The first formula follows by summing a single indicator function and applying Proposition~\ref{bas:mean}.
The forward direction of the second assertion follows from~\eqref{bas:counts} and $S_g(N-1)=o(N)$;
the converse follows immediately from $C_g(N)=N-1-4R_{-+}(N)$.
\end{proof}

\section{Proofs of the quantitative lower bounds and almost-all asymptotics}\label{sec:qplus}\label{sec:almost-all}

\subsection{Four decompositions of an integer}
We use a $q+1$ construction. Fix an integer $N\ge2$ and an odd prime $q$, and assume
\begin{equation}\label{plus:hypothesis}
q\nmid N,\qquad
p\mid q+1,\ p\text{ an odd prime}\Longrightarrow p\mid N.
\end{equation}
Write $b=q+1$ and $\beta=qb/2$. Since $b$ is even, $\beta$ is an integer.

\begin{lemma}\label{plus:four}
If
\begin{equation}\label{plus:A}
A\equiv N\pmod{q^2},\qquad 0<A<\frac N{q+1},
\end{equation}
define
\[
u=\frac{N-A}{q^2},\qquad v=N-\beta u,\qquad w=N-2v.
\]
Then $u,v,w$ are positive integers, and
\begin{equation}\label{plus:decompositions}
N=A+q^2u=\beta u+v=2v+w=bA+qw.
\end{equation}
\end{lemma}
\begin{proof}
The congruence and $0<A<N$ ensure that $u$ is a positive integer; $v,w$ are integers by definition. Substitution gives
\[
v=\frac{(q-1)N+bA}{2q}>0,\qquad
w=\frac{N-bA}{q}>0.
\]
The first three identities follow from the definitions, and the last from $bA+qw=N$. Every summand lies strictly in $(0,N)$.
\end{proof}

\begin{lemma}[The product of eight signs]\label{plus:negative}
Under~\eqref{plus:hypothesis} and~\eqref{plus:A},
\begin{equation}\label{plus:negative-product}
g(A)g(q^2u)g(\beta u)g(v)g(2v)g(w)g(bA)g(qw)=-1.
\end{equation}
\end{lemma}
\begin{proof}
For a prime $p$, let $e_p$ be the number of the eight displayed integers divisible by $p$, so that \(0\leq e_p\leq8\). The left-hand side equals $\prod_p(-1)^{e_p}$.

For $p=q$, we have $A\equiv N\not\equiv0\pmod q$, while $v=N-\beta u\equiv N\pmod q$ and $w=N-2v\equiv-N\pmod q$. Since $q\nmid2b$, exactly the three positions $q^2u,\beta u,qw$ are divisible by $q$. Hence $e_q=3$.

For $p=2$, each pair of summands has sum $N$. If $N$ is even, each pair has either zero or two even entries; if $N$ is odd, each has exactly one. Thus $e_2$ is always even.

If $p$ is an odd prime factor of $b$, our assumption gives $p\mid N$. Since $q\equiv-1\pmod p$ and $N=A+q^2u$, we have $A\equiv-u\pmod p$, so $p\mid A$ if and only if $p\mid u$. Also $p\mid\beta$, hence $p\mid v,w$. If $p\mid u$, all eight positions are divisible by $p$; otherwise, only $A,q^2u$ are not. Thus $e_p$ is $8$ or $6$, respectively.

Finally, if $p\nmid2qb$, pair the eight positions as
\[
(A,bA),\quad(q^2u,\beta u),\quad(v,2v),\quad(w,qw).
\]
The fixed coefficients in each pair are coprime to $p$, so divisibility occurs in pairs and $e_p$ is even.

This covers all primes. The only odd value of $e_p$ occurs at $q$, proving~\eqref{plus:negative-product}.
\end{proof}

Write $c_1,c_2,c_3,c_4\in\sgnset$ for the four correlation products. Lemma~\ref{plus:negative} shows that their product is $-1$, so both signs occur. This does not yet determine which same-sign pattern corresponds to a positive correlation product; that distinction will follow after using the summatory estimate.

Define
\[
K_q(N)=\#\{A\in\mathbb{N}:A\equiv N\pmod{q^2},\ 0<A<N/(q+1)\}.
\]
\begin{lemma}\label{plus:parameter-count}
Let $r$ be the least positive residue of $N$ modulo $q^2$, and let $B=\lfloor(N-1)/(q+1)\rfloor$. Then
\begin{equation}\label{plus:exact-K}
K_q(N)=\max\left(0,1+\left\lfloor\frac{B-r}{q^2}\right\rfloor\right),
\end{equation}
and
\begin{equation}\label{plus:lower-K}
K_q(N)\ge\frac{N}{(q+1)q^2}-1\ge\frac{N}{2q^3}-1.
\end{equation}
If $N\ge4q^3$, then $K_q(N)\ge N/(4q^3)$.
\end{lemma}
\begin{proof}
Since $q\nmid N$, we have $1\le r<q^2$. The congruence $A\equiv N\pmod{q^2}$ gives $A=r+kq^2$ for an integer $k\ge0$. The strict integer inequality $(q+1)A<N$ is equivalent to $A\le B$, so the parameters correspond exactly to integers $k$ satisfying
\[
0\le k\le\frac{B-r}{q^2}.
\]
There are no such integers if $B<r$; otherwise, there are $1+\lfloor(B-r)/q^2\rfloor$ of them. This proves~\eqref{plus:exact-K}.

Put $L=N/(q+1)$. If $L\le r$, then $K_q(N)=0\ge L/q^2-1$. If $L>r$, admissible $k$ satisfy $0\le k<(L-r)/q^2$, whence
\[
K_q(N)=\left\lceil\frac{L-r}{q^2}\right\rceil
\ge\frac{L-r}{q^2}>\frac L{q^2}-1.
\]
Now $q+1\le2q$ gives~\eqref{plus:lower-K}. Finally, if $N\ge4q^3$, then $N/(2q^3)-1\ge N/(4q^3)$.
\end{proof}

\begin{lemma}\label{plus:injection}
Denote the first summands in the four decompositions in~\eqref{plus:decompositions} by
\begin{align*}
a_1(A)&=A,& a_2(A)&=\frac{b(N-A)}{2q},\\
a_3(A)&=\frac{(q-1)N+bA}{q},&a_4(A)&=bA.
\end{align*}
Then the map $(A,i)\mapsto a_i(A)$ is injective on the following set:
\[
\{A\in\mathbb N:A\equiv N\pmod{q^2},\ 0<A<N/(q+1)\}
\times\{1,2,3,4\}.
\]

\end{lemma}
\begin{proof}
For each fixed $i$, the function $a_i(A)$ is affine with nonzero slope and is therefore injective. Using the integer expressions $a_2=\beta u$ and $a_3=2v$, we obtain
\[
a_1\equiv N,\quad a_2\equiv0,\quad
a_3\equiv2N,\quad a_4\equiv N\pmod q.
\]
Since $q$ is odd and $q\nmid N$, the residues $0,N,2N$ are distinct. It remains only to rule out equality between \(a_1\) and \(a_4\). For any $A,A'$ satisfying~\eqref{plus:A},
\[
a_1(A)\equiv N\pmod{q^2},\qquad
a_4(A')\equiv(q+1)N=N+qN\pmod{q^2}.
\]
Since $qN$ is not divisible by $q^2$, these cannot be equal. Thus all \(a_i(A)\) are distinct, and the corresponding ordered decompositions do not repeat.
\end{proof}

\begin{theorem}\label{plus:gap}
Under assumption~\eqref{plus:hypothesis},
\begin{equation}\label{plus:gap-exact}
N-1\pm C_g(N)\ge2K_q(N).
\end{equation}
If $N\ge4q^3$, then
\begin{equation}\label{plus:gap-size}
N-1\pm C_g(N)\ge\frac{N}{2q^3}.
\end{equation}
Each inequality holds separately for both choices of sign.
\end{theorem}
\begin{proof}
By Lemma~\ref{plus:negative}, each $A$ satisfying~\eqref{plus:A} yields at least one positive and one negative correlation product. Injectivity ensures that the decompositions from different \(A\) do not overlap, so $P_+(N),P_-(N)\ge K_q(N)$. Equation~\eqref{bas:product-counts} gives~\eqref{plus:gap-exact}, and Lemma~\ref{plus:parameter-count} gives~\eqref{plus:gap-size}.
\end{proof}

\begin{remark}
Assumption~\eqref{plus:hypothesis} is essential. For example, $q=5,N=76,A=1$ satisfy the parameter interval and congruence, but $3\mid q+1$ and $3\nmid N$. The eight positions are $(1,75,45,31,62,14,6,70)$, and the product of their $g$-values is $1$. Thus checking shared prime factors is a substantive part of the argument.
\end{remark}
In Subsection~\ref{sec:global}, we choose the least odd prime not dividing $N$ and show that it satisfies~\eqref{plus:hypothesis}. This makes the gap above applicable to all sufficiently large integers. We first turn to exponential sums to prove Theorem~\ref{thm:main-average}.

\subsection{Exponential sums and mean-square estimates}
For real $x\ge1$ and $\alpha\in\R$, write
\[
 G_x(\alpha)=\sum_{n\le x}g(n)\e(\alpha n),
 \qquad \e(y)=\exp(2\pi i y).
\]

Davenport\textquotesingle s classical theorem states that, for every fixed $A>0$,
\begin{equation}\label{avg:davenport}
 \sup_{\alpha\in\R}
 \left|\sum_{n\le x}\mu(n)\e(\alpha n)\right|
 \ll_A \frac{x}{(\log x)^A}
 \qquad(x\ge3).
\end{equation}
For the original paper, see~\cite{Davenport}; the precise form above appears in~\cite[p.~542, (1.1)]{GT}.

\begin{proposition}\label{avg:exponential}
For every fixed $A>0$,
\begin{equation}\label{avg:g-exponential}
 \sup_{\alpha\in\R}|G_x(\alpha)|
 \ll_A\frac{x}{(\log x)^A}
 \qquad(x\ge3).
\end{equation}
\end{proposition}

\begin{proof}
Using the convolution representation~\eqref{bas:convolution} and interchanging sums, we obtain
\[
 G_x(\alpha)
 =\sum_{d\le x}\kappa(d)
   \sum_{m\le x/d}\mu(m)\e(d\alpha m).
\]
We may assume $x\ge9$. Split the outer sum at $d=\sqrt x$.
For $d\le\sqrt x$, we have $x/d\ge3$ and
$\log(x/d)\ge\tfrac12\log x$, so~\eqref{avg:davenport} gives
\begin{align*}
 \left|\sum_{d\le\sqrt x}\kappa(d)
       \sum_{m\le x/d}\mu(m)\e(d\alpha m)\right|
 &\ll_A\frac{x}{(\log x)^A}
       \sum_{d\le\sqrt x}\frac{|\kappa(d)|}{d}\\
 &\ll_A H_1\frac{x}{(\log x)^A}.
\end{align*}
Here $H_1<\infty$ by~\eqref{bas:absolute}.
For the remaining part, the bound $|\mu(m)|\le1$ gives
\begin{align*}
 \left|\sum_{\sqrt x<d\le x}\kappa(d)
       \sum_{m\le x/d}\mu(m)\e(d\alpha m)\right|
 &\le x\sum_{d>\sqrt x}\frac{|\kappa(d)|}{d}\\
 &\le x^{7/8}
       \sum_{d\ge1}\frac{|\kappa(d)|}{d^{3/4}}
 =H_{3/4}x^{7/8}.
\end{align*}
For fixed $A$, the last term can be absorbed into
$O_A(x/(\log x)^A)$. All estimates are uniform in $\alpha$.
\end{proof}

\begin{proposition}\label{avg:mean-square}
For every fixed $A>0$,
\begin{equation}\label{avg:normalized-square}
 \sum_{2\le N\le X}\frac{|C_g(N)|^2}{N^2}
 \ll_A\frac{X}{(\log X)^A}
 \qquad(X\ge3).
\end{equation}
\end{proposition}

\begin{proof}
For an integer $T\ge2$, define
\[
 B_T(m)=\sum_{\substack{a+b=m\\1\le a,b\le T}}g(a)g(b).
\]
The $m$th Fourier coefficient of the polynomial $G_T(\alpha)^2$ is $B_T(m)$,
whose support is contained in $2\le m\le2T$. For $2\le m\le T+1$,
every decomposition $a+b=m$ into positive integers automatically has $a,b\le T$, so
\[
 B_T(m)=C_g(m)\qquad(2\le m\le T+1).
\]
Parseval\textquotesingle s identity gives
\[
 \sum_{m=2}^{2T}|B_T(m)|^2
 =\int_0^1|G_T(\alpha)|^4\,d\alpha.
\]
Since $|g(n)|=1$, another application of Parseval gives
\[
 \int_0^1|G_T(\alpha)|^2\,d\alpha=T.
\]
Applying Proposition~\ref{avg:exponential} with exponent $A/2$, we obtain
\begin{align}
 \sum_{m=2}^{T+1}|C_g(m)|^2
 &\le\int_0^1|G_T(\alpha)|^4\,d\alpha\notag\\
 &\le\|G_T\|_\infty^2
       \int_0^1|G_T(\alpha)|^2\,d\alpha
 \ll_A\frac{T^3}{(\log T)^A}.
 \label{avg:unweighted-square}
\end{align}
Thus, writing
\[
 A_0(u)=\sum_{2\le N\le u}|C_g(N)|^2,
\]
we have $A_0(u)\ll_A u^3/(\log u)^A$ for $u\ge3$.
Partial summation gives
\[
 \sum_{2\le N\le X}\frac{|C_g(N)|^2}{N^2}
 =\frac{A_0(X)}{X^2}
  +2\int_2^X\frac{A_0(u)}{u^3}\,du.
\]
Split the integral at $u=\sqrt X$. Use the trivial bound
$A_0(u)\ll u^3$ on the first part and~\eqref{avg:unweighted-square} on the second to obtain
\[
 \sum_{2\le N\le X}\frac{|C_g(N)|^2}{N^2}
 \ll_A\sqrt X+\frac{X}{(\log X)^A}
 \ll_A\frac{X}{(\log X)^A}.
\]
This proves the assertion.
\end{proof}

\begin{proof}[Proof of Theorem~\ref{thm:main-average}]
Let
\[
 \mathcal E_B(X)=
 \left\{2\le N\le X:
 |C_g(N)|>\frac{N}{(\log N)^B}\right\}.
\]
If $N\in\mathcal E_B(X)$, then
\[
 1<(\log X)^{2B}\frac{|C_g(N)|^2}{N^2}.
\]
Apply Proposition~\ref{avg:mean-square} with logarithmic saving exponent
$A+2B$ to obtain
\[
 |\mathcal E_B(X)|
 \le(\log X)^{2B}
     \sum_{2\le N\le X}\frac{|C_g(N)|^2}{N^2}
 \ll_{A,B}\frac{X}{(\log X)^A}.
\]
For integers outside this set, the pattern-counting identity~\eqref{bas:counts}, together with
the arbitrary fixed logarithmic saving from Proposition~\ref{bas:mean},
\[
 S_g(x)\ll_B\frac{x}{(\log x)^B},
\]
gives the result. The constant error from replacing $N-1$ by $N$
is absorbed into the stated error term.
\end{proof}


\subsection{The least missing odd prime and results for all large integers}\label{sec:global}
For $N\ge2$, define
\begin{equation}\label{glob:qdef}
q(N)=\min\{q:q\text{ is an odd prime and }q\nmid N\}.
\end{equation}
Existence follows directly by taking any prime factor of $2N+1$: it is odd and does not divide $N$.

We use the unconditional Rosser--Schoenfeld estimate~\cite[p.~70, (3.16)]{RS} for the Chebyshev function:
\begin{equation}\label{glob:theta}
\theta(x):=\sum_{p\le x}\log p
>x\left(1-\frac1{\log x}\right)\qquad(x\ge41).
\end{equation}

\begin{lemma}\label{glob:qbound}
Put $q=q(N)$. Then~\eqref{plus:hypothesis} holds. If $\log N\ge32$, then
\begin{equation}\label{glob:qsize}
q\le2\log N,\qquad N\ge4q^3.
\end{equation}
\end{lemma}
\begin{proof}
If an odd prime $p$ divides $q+1$, then $p\le(q+1)/2<q$. Minimality of $q$ implies $p\mid N$, so~\eqref{plus:hypothesis} holds.

Put $y=\log N\ge32$. If $q>2y$, every odd prime at most $2y$ divides $N$, whence
\[
\log N\ge\theta(2y)-\log2.
\]
Since $2y\ge64>41$ and $\log(2y)>4$, equation~\eqref{glob:theta} gives $\theta(2y)>3y/2$. This implies $y>3y/2-\log2$, contradicting $y\ge32$. Hence $q\le2y$.

The derivative of $e^y/y^3$ is $e^y(y-3)/y^4$, so this function is increasing for $y>3$. Also $e^{32}>2^{32}>32^4$; thus, for $y\ge32$,
\[
\frac{e^y}{y^3}\ge\frac{e^{32}}{32^3}>32.
\]
Together with $q\le2y$, this gives $N=e^y>32y^3\ge4q^3$.
\end{proof}

Applying Lemma~\ref{glob:qbound} to Theorem~\ref{plus:gap}, we obtain, for $\log N\ge32$,
\begin{equation}\label{glob:loggap}
N-1\pm C_g(N)\ge\frac{N}{16(\log N)^3}.
\end{equation}

\begin{proof}[Proof of Theorem~\ref{thm:main-log}]
By Proposition~\ref{bas:mean}, there is an absolute constant $A_4>0$ such that
\[
|S_g(x)|\le A_4\frac{x}{(\log x)^4}\qquad(x\ge3).
\]
Put $B_4=16A_4$. Since $\log(N-1)\ge\frac12\log N$ for $N\ge4$, we have
\begin{equation}\label{glob:Sbound}
|S_g(N-1)|\le B_4\frac{N}{(\log N)^4}\qquad(N\ge4).
\end{equation}
Define
\begin{equation}\label{glob:N1}
N_1=\left\lceil\exp\bigl(\max\{32,64B_4\}\bigr)\right\rceil.
\end{equation}
For $N\ge N_1$, equation~\eqref{glob:Sbound} gives
\[
2|S_g(N-1)|\le\frac{N}{32(\log N)^3}.
\]
By the pattern identities and~\eqref{glob:loggap},
\[
R_{--}(N),R_{++}(N)
\ge\frac14\left(\frac{N}{16(\log N)^3}-\frac{N}{32(\log N)^3}\right)
=\frac{N}{128(\log N)^3}.
\]
For the mixed-sign patterns, we directly have
\[
R_{-+}(N)=R_{+-}(N)
=\frac{N-1-C_g(N)}4
\ge\frac{N}{64(\log N)^3}.
\]
This completes the proof.
\end{proof}

Here $N_1$ is a threshold defined in terms of constants in classical mean-value estimates; we do not evaluate it numerically. Section~\ref{sec:explicit} proves existence separately, first treating $N\ge100000$ and then verifying the finite range below it, to obtain the optimal eventual threshold $18$. The number $100000$ separates the two parts of that argument, while $18$ is the existence threshold. Neither is established as a starting point for the quantitative lower bounds above.

\begin{corollary}[Linear lower bounds when the least missing prime is bounded]\label{glob:bounded-q}
Fix a real number $B\ge3$. There is a constant $N_B$ such that every $N\ge N_B$ with $q(N)\le B$ satisfies
\[
R_{--}(N),R_{++}(N)\ge\frac{N}{16B^3},\qquad
R_{-+}(N)=R_{+-}(N)\ge\frac{N}{8B^3}.
\]
\end{corollary}
\begin{proof}
We may take
\[
N_B=\left\lceil\max\left\{4B^3,
\exp\bigl((8B_4B^3)^{1/4}\bigr)\right\}\right\rceil.
\]
Theorem~\ref{plus:gap} gives $N-1\pm C_g(N)\ge N/(2B^3)$, while~\eqref{glob:Sbound} gives $2|S_g(N-1)|\le N/(4B^3)$. Substitution into~\eqref{bas:counts} proves the result.
\end{proof}

\begin{corollary}\label{glob:log-order}
For all $s,t\in\sgnset$, as $N\to\infty$ through all integers,
\[
R_{st}(N)=N^{1+o(1)}.
\]
\end{corollary}
\begin{proof}
By Theorem~\ref{thm:main-log} and $R_{st}(N)\le N-1$, we have
\[
1-\frac{3\log\log N+\log128}{\log N}
\le\frac{\log R_{st}(N)}{\log N}\le1.
\]
Both bounds tend to $1$.
\end{proof}
This determines the order on a logarithmic scale; it does not assert that $R_{st}(N)\sim N$.


\section{Goldbach-type problems for the Möbius function}\label{sec:mobius}
The ordinary Möbius function takes values in $\{-1,0,1\}$, so additive representations require simultaneous control of squarefreeness and signs.
We give exact expansions for the nine patterns and identify the part that remains unproved.
Counting binary additive representations by squarefree integers is a classical problem. For early work, see
Estermann~\cite{Estermann1931}; for exponential-sum methods and subsequent background, see Tolev~\cite{Tolev2005}.
\begin{definition}
For $N\ge2$ and $s,t\in\{-1,0,1\}$, write
\[
 B_{s,t}(N)=\sum_{a=1}^{N-1}\one_{\mu(a)=s}\one_{\mu(N-a)=t}.
\]
Define $A_1=6/\pi^2$ and
\[
 \CB(N)=\prod_p\left(1-\frac{\rho_N(p)}{p^2}\right),\qquad
 \rho_N(p)=\begin{cases}1,&p^2\mid N,\\2,&p^2\nmid N.\end{cases}
\]
\end{definition}
Modulo $p^2$, the forbidden residue classes are $a\equiv0,N$; these coincide when $p^2\mid N$.
The condition here is therefore $p^2\mid N$ and cannot be replaced by $p\mid N$.
The product converges absolutely, is positive, and satisfies
\[
 \CB(N)=\prod_p(1-2/p^2)
       \prod_{p^2\mid N}\frac{p^2-1}{p^2-2}.
\]

\subsection{Exact expansions for the nine patterns}
All sums below are over $1\le a<N$. Write
\begin{align*}
 M_0(N)&=N-1,&M_1(N)&=\sum_{a<N}\mu^2(a),\\
 M_2(N)&=\sum_{a<N}\mu^2(a)\mu^2(N-a),&S(N)&=\sum_{a<N}\mu(a),\\
 X_c(N)&=\sum_{a<N}\mu(a)\mu^2(N-a),&T(N)&=\sum_{a<N}\mu(a)\mu(N-a).
\end{align*}
Using $\one_{\mu(n)=0}=1-\mu^2(n)$ and
$\one_{\mu(n)=s}=(\mu^2(n)+s\mu(n))/2$ for $s=\pm1$, expand and apply $a\mapsto N-a$ to obtain
\begin{align}
 B_{0,0}&=M_0-2M_1+M_2,\label{eq:B00-expand}\\
 B_{0,t}=B_{t,0}&=\tfrac12(M_1-M_2)+\tfrac t2(S-X_c)\quad(t=\pm1),\label{eq:B0t-expand}\\
 B_{s,t}&=\tfrac14M_2+\tfrac{s+t}{4}X_c+\tfrac{st}{4}T\quad(s,t=\pm1).
 \label{eq:Bst-expand}
\end{align}
These identities also show directly that the nine actual counts sum to exactly $N-1$.

\subsection{Unconditional estimates}
\begin{proposition}\label{mob:known}
For every fixed $A>0$ and $\varepsilon>0$,
\begin{align*}
 M_1(N)&=A_1N+o(\sqrt N),\\
 M_2(N)&=\CB(N)N+O_\varepsilon(N^{2/3+\varepsilon}),\\
 S(N),X_c(N)&\ll_A N/(\log N)^A.
\end{align*}
The implied constants are independent of the integer $N$.
\end{proposition}
\begin{proof}
Summing the identity $\mu^2(n)=\sum_{d^2\mid n}\mu(d)$ first gives
$M_1(N)=\sum_{d\le\sqrt{N-1}}\mu(d)\lfloor(N-1)/d^2\rfloor$.
The floor-function error is $O(\sqrt N)$, and extending the main-term series to infinity also costs $O(\sqrt N)$,
while $\sum\mu(d)/d^2=1/\zeta(2)=A_1$.
To improve the error to $o(\sqrt N)$, we additionally use $M_\mu(x)=o(x)$.
Write $x=N-1$, fix $K>1$, and put $D=\lfloor\sqrt x/K\rfloor$.
The floor-function error for $d\le D$ is $O(D)$. Partial summation gives
$\sum_{d>D}\mu(d)/d^2=o(1/D)$, so extending the main term costs $o_K(\sqrt x)$.
The remaining actual sum can be written as
\[
 \sum_{d>D}\mu(d)\lfloor x/d^2\rfloor
 =\sum_{m\le x/(D+1)^2}
 \bigl(M_\mu(\sqrt{x/m})-M_\mu(D)\bigr)=o_K(\sqrt x),
\]
because, for fixed $K$, the outer sum has a bounded number of terms. Thus the limsup of the error divided by $\sqrt x$ is
$O(1/K)$. Letting $K\to\infty$ gives $o(\sqrt x)$.
Replacing $x=N-1$ by $N$ costs only $O(1)$.

To prove the second estimate, put $y=N^{1/3}$ and $Q_y(n)=\sum_{d\le y,d^2\mid n}\mu(d)$.
Using the elementary divisor bound $\tau(n)\ll_\eta n^\eta$ for any fixed $\eta>0$ (see~\cite{Apostol}), we have
\[
 \sum_{n<N}|\mu^2(n)-Q_y(n)|
 \le\sum_{y<d\le\sqrt N}(N/d^2+1)\ll N/y+\sqrt N.
\]
Also $|Q_y(n)|\le\tau(n)$. Hence truncating both squarefree indicators in $M_2$
introduces a total error of $O_\eta(N^\eta(N/y+\sqrt N))$.
The Chinese remainder theorem gives
\[
 \sum_{a=1}^{N-1}Q_y(a)Q_y(N-a)
 =N\!\!\sum_{\substack{d,e\le y\\(d,e)^2\mid N}}
 \frac{\mu(d)\mu(e)(d,e)^2}{d^2e^2}+O(y^2).
\]
The endpoint $a=N-1$ is retained here. To estimate the tail of the main-term series, write
$d=ku,e=kv$, where $k=(d,e)$ and $k^2\mid N$.
Dropping the coprimality conditions and Möbius factors and using $\sum_{u>z}u^{-2}\ll\min(1,z^{-1})$, we obtain
\[
 \sum_{\substack{d>y\ {\rm or}\ e>y\\(d,e)^2\mid N}}
 \frac{(d,e)^2}{d^2e^2}
 \ll \sum_{k^2\mid N}\frac1{k^2}\min(1,k/y)
 \ll\frac{\tau(N)}y.
\]
The double series converges absolutely. Its local factor at a prime, obtained from the exponent pairs $(0,0),(1,0),(0,1),(1,1)$, is
$1-2/p^2+\one_{p^2\mid N}/p^2$, so the infinite series equals $\CB(N)$.
Taking $\eta$ sufficiently small and substituting $y=N^{1/3}$ yields the stated error.

The bound for $S$ follows from the classical Möbius mean-value estimate. For $X_c$, expand the second factor and put
$D=(\log N)^{A+2}$. For $d\le D$, orthogonality of additive characters and
Davenport\textquotesingle s estimate~\eqref{avg:davenport} give the uniform bound
\[
 \left|\sum_{\substack{a<N\\a\equiv N\pmod{d^2}}}\mu(a)\right|
 \le\sup_\alpha\left|\sum_{a<N}\mu(a)\e(\alpha a)\right|
 \ll_L N/(\log N)^L.
\]
For $d>D$, use absolute values; the contribution is at most
$\sum_{D<d\le\sqrt N}(N/d^2+1)\ll N/D+\sqrt N$.
Taking $L>2A+3$ and combining the two parts yields $X_c(N)\ll_A N/(\log N)^A$.
The finitely many small values of $N$ can be absorbed into the implied constant.
\end{proof}
The error term obtained here for the binary squarefree count is
$O_\varepsilon(N^{2/3+\varepsilon})$, which suffices for the subsequent consequences in this section.
This error comes from estimating both the floor-function errors after truncation and the tail of large square factors. One cannot simply sum the floor-function errors in the untruncated double sum.

\subsection{The remaining pointwise correlation conjecture}
\begin{conjecture}\label{conj:unified}
As $N\to\infty$ through all integers, $T(N)=o(N)$.
\end{conjecture}
By Proposition~\ref{mob:known} and the three exact expansions, this conjecture is equivalent to the simultaneous asymptotics
\[
 B_{s,t}(N)=\kappa_{s,t}(N)N+o(N),\qquad
 \kappa_{s,t}(N)=\begin{cases}
 1-2A_1+\CB(N),&s=t=0,\\
 (A_1-\CB(N))/2,&\text{exactly one index is zero},\\
 \CB(N)/4,&s,t\ne0.
 \end{cases}
\]
The reverse implication follows immediately from $T=B_{+,+}+B_{-,-}-B_{+,-}-B_{-,+}$, since the four main terms cancel exactly.
These coefficients are nonnegative and sum to $1$; they also admit an interpretation through local exclusion probabilities modulo each $p^2$.
In fact, the five patterns with at least one zero index can already be handled unconditionally.
More precisely, Proposition~\ref{mob:known} gives, for any fixed $A>0$ and $\varepsilon>0$,
\begin{align*}
 B_{0,0}(N)&=(1-2A_1+\CB(N))N+O_\varepsilon(N^{2/3+\varepsilon}),\\
 B_{0,t}(N)=B_{t,0}(N)&=\tfrac12(A_1-\CB(N))N
                  +O_A(N/(\log N)^A)\quad(t=\pm1),\\
 B_{s,t}(N)&=\tfrac14\CB(N)N+\tfrac{st}{4}T(N)
                  +O_A(N/(\log N)^A)\quad(s,t=\pm1).
\end{align*}
To prove the last two formulas, first fix the exponent in the squarefree counting estimate at any value less than $1$,
and then absorb the resulting power error into any fixed logarithmic saving.
Thus the remaining quantity to control is $T(N)$, shared by the four nonzero patterns, rather than all nine terms.
If the remainders for all nine patterns were simultaneously improved to $O_\varepsilon(N^{1/2+\varepsilon})$,
then differences of the marginal counts would imply $S(N)=O_\varepsilon(N^{1/2+\varepsilon})$, and hence the Riemann hypothesis.
This formulation is stronger than Conjecture~\ref{conj:unified}.

\subsection{Relation to Chowla\textquotesingle s conjecture and the Corrádi--Kátai convolution conjecture}\label{sec:chowla}
Conjecture~\ref{conj:unified} has the random-sign motivation of a two-point Chowla-type conjecture, but fixed-shift correlations must be distinguished from binary additive convolution.
The classical source for Chowla\textquotesingle s conjecture is~\cite{Chowla1965}. Its standard fixed-shift two-point Möbius form is
\begin{equation}\label{mob:fixed-chowla}
 \sum_{n\le x}\mu(n)\mu(n+h)=o_h(x)
 \qquad(x\to\infty),\quad h\ge1\text{ fixed}.
\end{equation}
There is a corresponding formulation for the Liouville function. More generally, the fixed-affine-form version concerns
$\sum_{n\le x}\mu(a_1n+b_1)\mu(a_2n+b_2)$, with fixed coefficients,
$a_1,a_2>0$, and $a_1b_2-a_2b_1\ne0$; see~\cite{TaoChowla2016}.
Only the two-point form is stated here; the full higher-order Chowla conjecture is not limited to this assertion.

By contrast, our sum
\[
 T(N)=\sum_{a=1}^{N-1}\mu(a)\mu(N-a)
\]
uses the affine forms $a$ and $N-a$. The latter has slope $-1$, and its intercept and the summation interval both vary with $N$.
Equation~\eqref{mob:fixed-chowla} asserts a limit for each fixed $h$
and does not provide a uniform estimate that directly handles this varying family.
Consequently, equivalence between Conjecture~\ref{conj:unified} and the standard fixed-shift Chowla conjecture cannot be deduced merely by a change of variables.
In discussing the corresponding Liouville convolution, Mangerel~\cite[Section 1]{Mangerel2024} likewise emphasizes the difficulties caused by large shift parameters.

The distinction is transparent in Fourier analysis. For an integer $X\ge2$, write
$F_X(\alpha)=\sum_{n=1}^X\mu(n)\e(\alpha n)$. Orthogonality gives, respectively,
\begin{align*}
 T(N)&=\int_0^1 F_{N-1}(\alpha)^2\e(-N\alpha)\,d\alpha,\\
 \sum_{n=1}^{X-h}\mu(n)\mu(n+h)
 &=\int_0^1|F_X(\alpha)|^2\e(h\alpha)\,d\alpha\qquad(1\le h<X).
\end{align*}
The first identity extracts the sum of two indices; the second extracts their difference.
Both involve sign cancellation for multiplicative functions, but they address different coefficient problems.

For binary additive convolution, a more direct historical reference is the Corrádi--Kátai conjecture~\cite{CorradiKatai1969}:
\[
 \sum_{a=1}^{N-1}\lambda(a)\lambda(N-a)=o(N).
\]
For this attribution and subsequent developments, see~\cite[Section 1]{Mangerel2024} and
\cite[Section 1.2, Conjecture 1]{Krishnamoorthy}.
For work directly addressing Möbius additive convolution and its weighted averages, see also~\cite[Sections 1 and 5]{CGZ}.

\subsection{Average results and the pointwise requirement of our conjecture}\label{sec:chowla-average}
The different types of averaging in the literature must be distinguished.
Matomäki--Radziwiłł--Tao~\cite{MRTChowla2015} established Chowla-type results averaged over shift parameters;
Tao~\cite[Theorems 2 and 3]{TaoChowla2016} established logarithmically weighted two-point results for fixed affine forms.
For example, for every fixed $h\ge1$, the corresponding Möbius consequence is
\[
 \sum_{n\le x}\frac{\mu(n)\mu(n+h)}n=o_h(\log x).
\]
The weight $1/n$ cannot simply be removed, and the fixed $h$ cannot simply be replaced by a reflection parameter growing with the interval.
Tao--Teräväinen~\cite{TTScales2019} also proved certain unweighted fixed-shift correlations
outside a set of scales of logarithmic density zero. This still does not assert the conclusion for every sufficiently large target integer, as required here.

For our $T(N)$, Davenport\textquotesingle s estimate directly yields an average result in the \emph{target sum $N$}.
This differs from a fixed-shift assertion, averaging over shifts, or averaging over logarithmic scales.
\begin{proposition}[A mean-square estimate for the reflected correlation]\label{mob:mean-square}
For any fixed $L>0$,
\[
 \sum_{2\le N\le X}\frac{|T(N)|^2}{N^2}
 \ll_L\frac{X}{(\log X)^L}\qquad(X\ge3).
\]
Consequently, for any fixed $A,B>0$, apart from at most $O_{A,B}(X/(\log X)^A)$ integers
$2\le N\le X$, the nine patterns simultaneously satisfy
\[
 B_{s,t}(N)=\kappa_{s,t}(N)N+
 O_B\left(\frac{N}{(\log N)^B}\right),\qquad s,t\in\{-1,0,1\}.
\]
\end{proposition}
\begin{proof}
By Davenport\textquotesingle s estimate, for any fixed $L>0$,
$\|F_X\|_\infty\ll_L X/(\log X)^{L/2}$.
On the other hand, Parseval gives $\int_0^1|F_X(\alpha)|^2\,d\alpha=\sum_{n\le X}\mu^2(n)\le X$.
The coefficient of the squared polynomial $F_X^2$ at frequency $2\le N\le X+1$ is exactly $T(N)$. Hence
\[
 \sum_{2\le N\le X+1}|T(N)|^2
 \le\int_0^1|F_X(\alpha)|^4\,d\alpha
 \le\|F_X\|_\infty^2\int_0^1|F_X(\alpha)|^2\,d\alpha
 \ll_L\frac{X^3}{(\log X)^L}.
\]
The same partial-summation argument as in Proposition~\ref{avg:mean-square} gives the normalized mean-square bound.
Apply this bound to integers satisfying $|T(N)|>N/(\log N)^B$, taking exponent $A+2B$,
to obtain the stated exceptional-set estimate. Finally, substitute into the nine-pattern expansions and apply Proposition~\ref{mob:known}.
\end{proof}
This average consequence establishes our conjecture for almost all target integers, while still allowing infinitely many sparse exceptions.
It does not replace the assertion $T(N)=o(N)$ as $N\to\infty$ through all integers.


\section{Proof of explicit existence and the complete exceptional sets}\label{sec:explicit}

We prove Theorem~\ref{thm:main-explicit} using a new approach. Retaining the notation
$A_s=\{a\in\mathbb N:g(a)=s\}$ from the introduction, we determine which integers belong to all four sumsets simultaneously.
Equivalently, we determine when
\[
N\in\bigcap_{s,t\in\sgnset}(A_s+A_t).
\]
For $N\ge100000$, we prove that this membership always holds. For $2\le N<100000$,
we verify each integer to determine all exceptions listed in the introduction. The dividing point used here is $100000$; changing it requires rechecking the corresponding covering chains, witnesses, and strict endpoints, rather than merely substituting a different number.
The number $100000$ serves only to divide the two parts of the argument. Both parts require computer-assisted verification.

The infinite range is classified according to $q=q(N)$, defined in Subsection~\ref{sec:global}.
Since every odd prime less than $q$ divides $N$, larger $q$ forces $N$ to contain more small prime factors.
We treat, in order: $q=3,5,7$; $11\le q\le59$ once $N$ reaches the corresponding threshold; and all $q\ge61$.
The integers not yet covered form an explicitly enumerable finite set, and a concrete decomposition is checked for each of them.
Thus the task is to cover the infinite range by mathematical arguments and enumerate all remaining integers exactly, rather than extrapolate from computation up to a finite bound.

We first reduce additive representation to finding parameters of both signs, and then explain why covering chains require only finite data.
All verification conditions and data versions are given in Appendices~\ref{app:certificates} and~\ref{app:small-range}.

\subsection{Four integer decompositions}

For an integer $N\ge2$ and an odd prime $q$, write
\[
 t=q-1,\qquad d=\rad(t).
\]
Define the binary relation $H(N,q)$ by
\begin{equation}\label{exp:H}
 H(N,q)\quad\Longleftrightarrow\quad
 \begin{cases}
 q\nmid N,\\
 p\mid q-1\Longrightarrow p\mid N&\text{for every odd prime }p.
 \end{cases}
\end{equation}
Here $q$ is a variable parameter in the relation and is not prescribed to be the previously defined $q(N)$.
Equation~\eqref{glob:qdef} defines the least odd prime $q(N)$ not dividing $N$ and guarantees its existence.
By definition, $q(N)\nmid N$. If an odd prime $p$ divides $q(N)-1$, then $p<q(N)$, so minimality gives $p\mid N$.
Thus $H(N,q(N))$ holds. The results below stated for $N,q$ satisfying $H(N,q)$ allow any such $q$; we choose $q=q(N)$ only when its minimality is needed.

\begin{lemma}\label{exp:four}
Suppose that $N,q$ satisfy the relation $H(N,q)$. If
\begin{equation}\label{exp:Aconditions}
 A\equiv N\pmod{q^2},\qquad 0<A<N/t^2,
\end{equation}
then setting
\[
 u=(N-A)/q^2,\qquad v=N-tqu,\qquad w=N-tv
\]
yields four decompositions of \(N\) into positive integers:
\begin{equation}\label{exp:fouridentities}
 N=A+q^2u=tqu+v=tv+w=t^2A+qw.
\end{equation}
Moreover, $q\nmid v,w$, and
\[
 g(q^2u)=g(qu),\qquad g(qw)=-g(w).
\]
\end{lemma}
\begin{proof}
The congruence makes $u$ an integer, and $A<N/t^2\le N$ gives $u>0$.
Direct calculation gives
\[
 v=\frac{N+tA}{q}>0,\qquad
 w=\frac{N-t^2A}{q}>0.
\]
These integers satisfy~\eqref{exp:fouridentities}; the last identity can also be checked directly from
$t^2+q(1-t)=1$.
The definitions give $v\equiv N\pmod q$ and $w\equiv2N\pmod q$.
Since $q$ is odd and $H(N,q)$ gives $q\nmid N$, we have $q\nmid v,w$.
The integers $qu$ and $q^2u$ have the same distinct prime factors, even when $q\mid u$;
also $q\nmid w$. The final two sign identities follow.
\end{proof}

The four decompositions themselves imply the mixed-sign patterns, but do not yet guarantee both same-sign patterns. We will show that if one same-sign pattern is absent,
then all parameters satisfying specified conditions are forced to have the same sign.
Finding both a positive and a negative sign among these parameters therefore contradicts the absence assumption and establishes the same-sign patterns as well.

\Needspace{6\baselineskip}
\begin{proposition}\label{exp:rigidity}
Suppose that $N,q$ satisfy the relation $H(N,q)$, and let $\varepsilon\in\{-1,1\}$.
\begin{enumerate}
\item If $N$ is odd and the pattern $(\varepsilon,\varepsilon)$ is absent,
then every $A$ satisfying~\eqref{exp:Aconditions} and $d\mid A$ has
$g(A)=-\varepsilon$.
\item If $N$ is even, put $m=tqd$. If the pattern $(\varepsilon,\varepsilon)$ is absent,
then every integer $v$ satisfying
\begin{equation}\label{exp:vconditions}
 N/q<v<N/t,\qquad v\equiv N-m/2\pmod m
\end{equation}
has $g(v)=-\varepsilon$.
\end{enumerate}
In either case, the existence of one admissible parameter guarantees both ordered mixed-sign patterns.
\end{proposition}
\begin{proof}
First suppose $N$ is odd. Since $2\mid d$ and $d\mid A$, the integer $A$ is even,
so $u$ is odd; since $t$ is even, $v$ is odd as well.
If an odd prime $p$ divides $t$, the relation $H(N,q)$ and the condition $d\mid A$ give $p\mid N,A$, and $p\ne q$,
hence $p\mid u,v$. Thus, when $qu$ or $v$ is multiplied by $t$,
all odd prime factors of $t$ are already present, and the only newly introduced distinct prime factor is $2$;
multiplying $A$ by $t^2$ introduces no distinct prime factor.
Write
\[
 \alpha=g(A),\quad x=g(q^2u),\quad y=g(v),\quad z=g(w).
\]
The signs in the four decompositions are, in order,
\begin{equation}\label{exp:oddsigns}
 (\alpha,x),\quad(-x,y),\quad(-y,z),\quad(\alpha,-z).
\end{equation}
If $\alpha=\varepsilon$, the absence assumption successively forces
$x=-\varepsilon$, $y=-\varepsilon$, $z=-\varepsilon$, 
but the last pair is then $(\varepsilon,\varepsilon)$, a contradiction.

Now suppose $N$ is even and $v$ satisfies~\eqref{exp:vconditions}.
Put
\[
 u=\frac{N-v}{tq},\qquad A=N-q^2u=\frac{qv-N}{t}.
\]
The congruence ensures integrality and gives
\begin{equation}\label{exp:uclass}
 u\equiv d/2\pmod d.
\end{equation}
The inequalities $N/q<v<N/t$ imply $u>0$ and $0<A<N/t^2$.
Since $d$ is even and squarefree, $d/2$ is odd, so~\eqref{exp:uclass}
shows that $u$ is odd and contains all odd prime factors of $t$.
The relation $H(N,q)$, together with the evenness of $N$, gives $d\mid N$.
Also $d\mid m/2$, so $d\mid v$.
The integer $A=N-q^2u$ is odd and also contains all odd prime factors of $t$.
The four sign pairs are therefore
\begin{equation}\label{exp:evensigns}
 (\alpha,x),\quad(-x,y),\quad(y,z),\quad(-\alpha,-z).
\end{equation}
Here multiplying $qu$ by $t$ introduces only the prime $2$, multiplying $v$ by $t$ leaves its prime-factor set unchanged,
and multiplying $A$ by $t^2$ introduces only the prime $2$.
If $y=\varepsilon$, the second, third, and first pairs successively force
$x=\varepsilon$, $z=-\varepsilon$, $\alpha=-\varepsilon$, 
and the fourth pair becomes $(\varepsilon,\varepsilon)$, a contradiction.

In both sign tables, the product of the four pairwise products is $-1$.
Thus at least one pair has product $-1$, giving a mixed-sign decomposition;
exchanging its summands gives the other ordered mixed-sign pattern.
\end{proof}

Therefore, for odd $N$, it suffices to find both signs among the admissible $A$;
for even $N$, it suffices to find both signs among the admissible $v$.
In the odd branch, $d\mid A$ is an additional parameter restriction used to control prime-factor changes. It does not follow automatically from
$A\equiv N\pmod{q^2}$. Since $(d,q^2)=1$, the two congruence conditions are compatible;
the existence of suitable parameters in the required short interval must still be guaranteed by the witnesses or counting arguments below.

\subsection{Covering an infinite range with finite chains}

Consider the even branch. A parameter must have the required sign, lie in the required residue class, and belong to $(N/q,N/t)$.
For a fixed integer $n$, the condition $N/q<n<N/t$ is equivalent to $tn<N<qn$.
Thus each candidate parameter $n$ gives an interval $(tn,qn)$ of target integers that it handles.
If the intervals arising from a sequence of parameters in the same residue class and with the same sign overlap consecutively, they cover a whole range of $N$.
The following lemma further shows that, when the parameters already contain a prime factor $\ell$, multiplication by powers of $\ell$ preserves their signs,
allowing finitely many chains to cover all larger integers.

\begin{lemma}[Scaling covering chains]\label{exp:chain}
Let $q\ge3$, put $t=q-1$, and let $T>0$. Let $m$ be a positive integer and $\ell$ a prime coprime to $m$.
Let $\mathcal R$ be a set of residue classes modulo $m$ invariant under multiplication by $\ell^{-1}$.
Suppose that for each $a\in\mathcal R$ and each sign $\varepsilon\in\{-1,1\}$ there is a finite chain of positive integers
\[
n_{a,\varepsilon,0}<\cdots<n_{a,\varepsilon,J(a,\varepsilon)},
\]
whose elements satisfy
\begin{gather}
n_{a,\varepsilon,i}\equiv a\pmod m,\qquad
\ell\mid n_{a,\varepsilon,i},\qquad
g(n_{a,\varepsilon,i})=\varepsilon,
\label{exp:chainsupport}\\
t n_{a,\varepsilon,0}<T<q n_{a,\varepsilon,0},\qquad
t n_{a,\varepsilon,i+1}<q n_{a,\varepsilon,i},\qquad
q n_{a,\varepsilon,J(a,\varepsilon)}>\ell T,
\label{exp:chainends}
\end{gather}
where the middle inequality holds for $0\le i<J(a,\varepsilon)$.
Then, for every real $X\ge T$, the interval $(X/q,X/t)$ contains integers of both signs in every class $r\in\mathcal R$.
\end{lemma}
\begin{proof}
Fix the target class $r\in\mathcal R$ and target sign $\varepsilon$. Put
\[
k=\left\lfloor\frac{\log(X/T)}{\log\ell}\right\rfloor\ge0,
\qquad Y=X/\ell^k\in[T,\ell T).
\]
Since $(\ell,m)=1$, the integer $\ell^k$ is invertible modulo $m$. Let
\begin{equation}\label{exp:chain-target-class}
a\equiv\ell^{-k}r\pmod m.
\end{equation}
Invariance of $\mathcal R$ implies $a\in\mathcal R$, so the hypothesis supplies a chain for $(a,\varepsilon)$.
Abbreviate this particular chain as $n_0<\cdots<n_J$. Every element satisfies
$n_i\equiv a\pmod m$ and $g(n_i)=\varepsilon$.

Since $n_i$ is strictly increasing and $tn_{i+1}<qn_i$, the left and right endpoints of the open intervals
$(tn_i,qn_i)$ both increase, and consecutive intervals have a nonempty open intersection.
Hence
\[
\bigcup_{i=0}^J(tn_i,qn_i)=(tn_0,qn_J)\supset[T,\ell T].
\]
The inclusion at both ends follows from the strict inequalities in~\eqref{exp:chainends}, so neither endpoint of the base interval is omitted.
Since $Y\in[T,\ell T)$, some element $n_i$ of the selected chain satisfies
\[
tn_i<Y<qn_i.
\]
Put $v=\ell^k n_i$. Then
\[
X/q<v<X/t,\qquad
v\equiv\ell^k a\equiv r\pmod m.
\]
Since $\ell\mid n_i$, multiplication by $\ell^k$ introduces no new distinct prime factor, so
\[
g(v)=g(n_i)=\varepsilon.
\]
This proves the assertion. The number $Y$ need not be an integer. The chain length $J(a,\varepsilon)$ is part of the fixed finite data;
it may depend on the residue class and sign, but need not increase with $X$.
\end{proof}

The covering-chain lemma specifies the exact conditions to check. For finitely many intermediate primes $q$,
we next choose an auxiliary prime $\ell$ and list chains satisfying these conditions.
A single finite verification then handles every integer above the threshold in each branch.

\begin{proposition}\label{exp:middle}
Let
\[
 \mathcal Q=\{11,13,17,19,23,29,31,37,41,43,47,53,59\}.
\]
For $q\in\mathcal Q$, suppose that $N,q$ satisfy the relation $H(N,q)$. The condition
\begin{equation}\label{exp:T}
 N\ge T(q):=128dq^2t^2
\end{equation}
guarantees the occurrence of all four sign patterns.
\end{proposition}
\begin{proof}[Proof using fixed finite verification data]
Put $h_0=tq$ and $m=dh_0$. For each $q\in\mathcal Q$, choose the auxiliary prime
$\ell$ listed in the parameter table below; direct checks give $(\ell,m)=1$. For every $r\bmod h_0$ with $q\nmid r$
and every $\varepsilon\in\{-1,1\}$, the finite verification data also supply a chain of positive integers $j_0<\cdots<j_J$
satisfying
\begin{gather}
 \ell j_i\equiv r\pmod{h_0},\qquad g(d\ell j_i)=\varepsilon,
 \label{exp:middleclass}\\
 td\ell j_0<T(q)<qd\ell j_0,\quad
 tj_{i+1}<qj_i,\quad qd\ell j_J>\ell T(q).
 \label{exp:middleends}
\end{gather}
Put $n_i=d\ell j_i$. Residue classes $r_0=dR\bmod m$ correspond bijectively to
$R\bmod h_0$; more precisely,
\[
 d\ell j\equiv dR\pmod{dh_0}
 \iff \ell j\equiv R\pmod{h_0}.
\]
Thus we do not need to assume that $d$ is invertible modulo $m$.
The classes covered are exactly
\[
 \mathcal R=\{r_0\bmod m:d\mid r_0,\ q\nmid r_0\}.
 \]
Since $(\ell,m)=1$ and $q\mid m$, this set is invariant under multiplication by $\ell^{-1}$.
Lemma~\ref{exp:chain} therefore applies.

If $N$ is even, then $d\mid N$ and $d\mid m/2$,
so the target class $N-m/2$ is divisible by $d$.
Also $q\mid m/2$ and $q\nmid N$, so this target class is not divisible by $q$.
Thus both signs occur among the admissible $v$, and Proposition~\ref{exp:rigidity} gives all patterns.

If $N$ is odd, the same finite verification data provide, for every $r\bmod q^2$ with $q\nmid r$
and every $\varepsilon$, an integer $k$ satisfying
\begin{equation}\label{exp:middleodd}
 0<k<128q^2,\qquad k\equiv r\pmod{q^2},\qquad g(dk)=\varepsilon.
\end{equation}
Since $(d,q)=1$, choose $r\equiv Nd^{-1}\pmod{q^2}$ and $A=dk$.
Then $d\mid A$, $A\equiv N\pmod{q^2}$, and
\[
 0<A<128dq^2=T(q)/t^2\le N/t^2.
\]
Proposition~\ref{exp:rigidity} again gives the result.

The finite verification data contain $1770208$ even-branch chain coefficients and $32492$ odd-branch witnesses.
Complete residue-class coverage, signs, coprimality conditions, and all strict endpoints have been checked individually.
Appendix~\ref{app:certificates} specifies the exact meaning of this finite verification;
these integer data form part of the proof of this proposition.
\end{proof}

\subsection{The cases with auxiliary prime \texorpdfstring{$3,5,7$}{3,5,7}}

We treat the three smallest branches. For $q=3$, we list the covering chains in full to illustrate the choice of residue classes, signs, and chain elements;
the cases $q=5,7$ use the same mechanism, with separate checks of their specific divisibility conditions.

When $q=3$ and $N$ is even, the modulus in~\eqref{exp:vconditions} is $12$,
and the target residue class $r\equiv N-6\pmod{12}$ belongs to
$\mathcal R=\{2,4,8,10\}$. All eight chains for this branch are listed below.
\begin{center}
\small
\begin{tabular}{rrl}
\toprule
$a\bmod12$&$\varepsilon$&The four integers in the chain\\
\midrule
2&$-1$&$49910,\ 74810,\ 112130,\ 168110$\\
2&$+1$&$49970,\ 74930,\ 112370,\ 168530$\\
4&$-1$&$49960,\ 74920,\ 112360,\ 168460$\\
4&$+1$&$49840,\ 74740,\ 112060,\ 167980$\\
8&$-1$&$49760,\ 74600,\ 111860,\ 167780$\\
8&$+1$&$49940,\ 74900,\ 112340,\ 168260$\\
10&$-1$&$49990,\ 74950,\ 112390,\ 168550$\\
10&$+1$&$49810,\ 74710,\ 112030,\ 168010$\\
\bottomrule
\end{tabular}
\end{center}
Factoring the entries individually verifies that every element of each chain is congruent to $a$ modulo $12$, is divisible by $5$,
and has the listed $g$-value. Writing each row as $n_0,n_1,n_2,n_3$, all eight chains satisfy
\[
2n_0<100000<3n_0,\qquad
2n_{i+1}<3n_i\ (0\le i<3),\qquad
3n_3>500000.
\]
The minimum strict margins over all chains and links are
\begin{align*}
\min(100000-2n_0)&=20,&\min(3n_0-100000)&=49280,\\
\min(3n_3-500000)&=3340,&\min(3n_i-2n_{i+1})&=10.
\end{align*}
Each minimum is taken over all relevant rows, and the last also over all adjacent positions.
Thus, for each chain, the union of the open intervals $(2n_i,3n_i)$ covers all of $[100000,500000]$.
Since $5^{-1}\equiv5\pmod{12}$, multiplication by $5^{-1}$ interchanges
$2\leftrightarrow10$ and $4\leftrightarrow8$, preserving the set of residue classes.
Applying Lemma~\ref{exp:chain} with $T=100000$ and $\ell=5$ guarantees that, for every even integer
$N\ge100000$ with $3\nmid N$, both signs occur among admissible $v$.
Proposition~\ref{exp:rigidity} then gives all four patterns.

For example, when $N=1000000$, the target class is $r=10\bmod12$. Take $k=1$; then
$Y=200000$, so before scaling we choose the chain in class $a\equiv5^{-1}r\equiv2\pmod{12}$.
In the negative-sign chain, $n_1=74810$ satisfies
\[
149620=2n_1<200000<3n_1=224430.
\]
Hence $v_-=5n_1=374050$ satisfies the required interval, congruence, and negative-sign conditions.
Likewise, the positive-sign chain supplies $74930$, giving $v_+=374650$.
The selection first fixes the entire chain corresponding to the residue class and sign, and then uses interval coverage to select one of its elements;
it does not arbitrarily change a chain from one fixed residue class to another.

When $q=5$ and $N$ is even, the modulus is $40$, and the required residue classes are all even classes not divisible by
$5$. The thirty-two chains contain 181 integers in total, all divisible by $3$,
and their intervals $(4n_i,5n_i)$ cover $[100000,300000]$.
Taking $\ell=3$ completes the scaling to the infinite range. Both families of residue classes are invariant under their respective $\ell^{-1}$.

The complete even witnesses modulo $9$ for the odd branch are listed in Appendix~\ref{cert:smalltables};
the largest is $64$, and $4\cdot64<100000$.
Every admissible class modulo $25$ likewise has even witnesses of both signs; the largest is $312$,
and $16\cdot312<100000$. This completes the branches where the least missing odd prime is $3$ or $5$.

Now let $q=7$ be the least missing odd prime. Then $15\mid N$ and $7\nmid N$,
so $N,7$ satisfy the relation $H(N,7)$. In the even branch, $t=d=6$, $m=252$, and the target conditions are
\[
 N/7<v<N/6,\qquad v\equiv N-126\pmod{252}.
\]
In particular,
\[
 u=(N-v)/42\equiv3\pmod6.
\]
Thus $3\mid u,A$ and $6\mid v$. Replacing $49u$ by $42u$ introduces only the prime $2$;
replacing $v$ by $6v$ leaves the distinct prime-factor set unchanged;
and replacing $A$ by $36A$ introduces only the prime $2$.
This directly rechecks that the shared factor $3$ does not invalidate~\eqref{exp:evensigns}.

For each $r\bmod252$ with $6\mid r$ and $7\nmid r$, and each sign,
the finite verification data provide chains covering $[100000,1000000]$ and $[1000000,5000000]$, respectively.
There are 144 chains and 1990 integers in total.
All integers in the second group are divisible by $5$, so $\ell=5$ extends coverage to all larger real targets;
the first group covers the smaller interval. Proposition~\ref{exp:rigidity} completes the even branch.

For odd $N$, choose $6\mid A$ and $A\equiv N\pmod{49}$.
All negative-sign witnesses in the fixed table are at most $1074$;
all positive-sign witnesses outside class $19\bmod49$ are at most $2777$.
Since $36\cdot2777=99972<100000$, these all satisfy the strict upper bound.
The positive-sign witness for class $19$ is $2910$, valid for $N>104760$.
The remaining conditions are
\[
 100000\le N\le104760,\quad 15\mid N,\quad N\ \text{odd},
 \quad N\equiv19\pmod{49}.
\]
These congruences are equivalent to $N\equiv705\pmod{1470}$;
the interval contains exactly the three integers $100665,102135,103605$.
They have positive-positive representations
\begin{align*}
 100665&=1+2^3\cdot12583,\\
 102135&=2\cdot3^2+3\cdot34039,\\
 103605&=2\cdot3+3^4\cdot1279.
\end{align*}
Here $12583,34039,1279$ are primes, as can be checked by trial division up to $112,184,35$, respectively.
Both summands in each identity have positive $g$-value.
An admissible negative-sign $A$ already guarantees the negative-negative pattern: if that pattern were absent,
Proposition~\ref{exp:rigidity} would force all admissible $A$ to be positive, a contradiction.
Any admissible $A$ also guarantees both mixed-sign patterns. Thus these three cases are fully handled.

\subsection{A uniform treatment of larger auxiliary primes}

It remains to handle arbitrarily large $q$. Adding auxiliary primes and covering chains one by one would not allow the proof to end with finite data.
Instead, we use a uniform interval-length condition: a sufficiently long interval contains both signs in the same residue class.
We then exploit the prime divisibility forced by $q=q(N)$ to show that this condition is automatic for all $q\ge61$.

\begin{lemma}[A three-point product identity]\label{exp:triple}
If $c,j,k$ are positive integers with $(j,k)=1$, then
\[
 g(cj^2)g(cjk)g(ck^2)=g(c).
\]
\end{lemma}
\begin{proof}
If a prime $p$ divides $c$, it occurs in all three positions and contributes $-1$.
If $p\nmid c$, the condition $(j,k)=1$ implies that it divides at most one of $j,k$;
if it occurs, it appears in exactly two positions and contributes $1$.
No condition $(c,jk)=1$ is required.
\end{proof}

\begin{lemma}[A three-point lemma allowing nonunit square congruences]\label{exp:square}
Let $c,m$ be positive integers and $0\le L<U$. Write $a=(c,m)$ and $M=m/a$.
If there is an integer $h$ satisfying
\[
 (h,M)=1,\qquad ch^2\equiv r\pmod m,
\]
and $\sqrt U-\sqrt L>2M\sqrt c$, then the class $r$ modulo $m$ contains an integer in $(L,U)$
with sign $g(c)$.
\end{lemma}
\begin{proof}
Choose the least integer $j\equiv h\pmod M$ strictly greater than $\sqrt{L/c}$,
and put $k=j+M$. Then
\[
 j\le\sqrt{L/c}+M,\qquad
 k\le\sqrt{L/c}+2M<\sqrt{U/c},\qquad (j,k)=(h,M)=1.
\]
Thus $cj^2,cjk,ck^2$ all lie strictly in $(L,U)$.
Since $cM=(c/a)m$, all three are congruent to $ch^2$ modulo $m$.
Lemma~\ref{exp:triple} shows that at least one has sign $g(c)$;
otherwise their sign product would be $-g(c)$, a contradiction.
\end{proof}

\begin{lemma}[Both signs in the same residue class]\label{exp:length}
Let $d,m$ be positive integers, let $r$ be an integer, and let $0\le L<U$. Suppose that $d\mid m$ and $d\mid r$, and let $p$ be a prime not dividing $m$.
If
\begin{equation}\label{exp:lengthcondition}
 \sqrt U-\sqrt L>2(m/d)p\sqrt{2m},
\end{equation}
then the class $r$ modulo $m$ contains integers of both signs in $(L,U)$.
\end{lemma}
\begin{proof}
Choose the least positive representative $1\le b\le m$ of the class.
Since $p\nmid m$, at least one of $b,b+m$ is not divisible by $p$.
Choose one of them as $c$ with $0<c<2m$; if $b=m$, simply take $c=m$.
Then $d\mid c$ and $p\nmid c$. Write $a=(c,m)$ and $M=m/a$.
For $c$, Lemma~\ref{exp:square} applies with $h=1$.
For $c'=cp^2$, we have $(c',m)=a$ and $g(c')=-g(c)$.
If $M>1$, choose $ph\equiv1\pmod M$; if $M=1$, take $h=1$.
Then $(h,M)=1$, and $cp^2h^2-c$ is divisible by $cM$ and hence by $m$.
Thus $c'h^2\equiv r\pmod m$.
The integer $c'$ itself need not be congruent to $r$.
The two applications of the three-point lemma require lengths $2M\sqrt c$ and $2Mp\sqrt c$, respectively;
both are strictly smaller than the right-hand side of~\eqref{exp:lengthcondition}
because $a\ge d$ and $c<2m$. Hence both opposite signs occur.
\end{proof}

\begin{proposition}[A general sufficient condition]\label{exp:B}
Suppose that $N,q$ satisfy the relation $H(N,q)$. If a prime $p$ satisfies $p\nmid q(q-1)$ and
\begin{equation}\label{exp:Bcondition}
 N>B(q,p):=32p^2t^4q^5d,
\end{equation}
then all four patterns occur.
\end{proposition}
\begin{proof}
For odd $N$, since $(d,q)=1$, the conditions
$d\mid A$, $A\equiv N\pmod{q^2}$
determine a unique class modulo $dq^2$, and this class is divisible by $d$.
Apply Lemma~\ref{exp:length} with
$m=dq^2$ and $(L,U)=(0,N/t^2)$.
The required condition is
\[
 N>8p^2dt^2q^6.
\]
Condition~\eqref{exp:Bcondition} is stronger, since $4t^2>q$.
For even $N$, put $m=tqd$. The target class $N-m/2$ is divisible by $d$,
and $p\nmid m$. Observe that
\[
 \sqrt{N/t}-\sqrt{N/q}
 =\frac{\sqrt N}{\sqrt{tq}(\sqrt q+\sqrt t)}
 >\frac{\sqrt N}{2q\sqrt t}.
\]
Condition~\eqref{exp:Bcondition} makes the right-hand side greater than
$2tqp\sqrt{2tqd}=2(m/d)p\sqrt{2m}$.
Both cases now follow from Lemma~\ref{exp:length} and Proposition~\ref{exp:rigidity}.
\end{proof}

Retain $q(N)$ from~\eqref{glob:qdef}, and write
\[
 P(q)=\prod_{\substack{3\le r<q\\r\ {\rm prime}}}r.
\]
When $q=q(N)$, we have $P(q)\mid N$, and $N,q$ satisfy $H(N,q)$.

Now return to $q=q(N)$. The sufficient condition above requires $N$ to exceed $B(q,p)$;
on the other hand, $P(q)\mid N$ gives $N\ge P(q)$.
Thus it suffices to compare $P(q)$ with $B(q,p)$.
Nagura\textquotesingle s theorem has two roles here: it supplies a moderately sized $p>q$ and controls the ratio of consecutive primes,
allowing inequalities at finitely many base points to be propagated along all subsequent primes.

\begin{proposition}[Completion of the infinite tail]\label{exp:tail}
If $q(N)\ge61$, then all four patterns occur.
\end{proposition}
\begin{proof}
For $q=61$ take $p=7$, and for $q=67$ take $p=5$. Both choices satisfy $p\nmid q(q-1)$.
Exact integer calculations give
\begin{align*}
 P(61)&=961380175077106319535\\
      &>514898337587558400000=B(61,7),\\
 P(67)&=58644190679703485491635\\
      &>1352644522537268505600=B(67,5).
\end{align*}
Thus $N\ge P(q)>B(q,p)$, and Proposition~\ref{exp:B} applies.

For all $q\ge71$, use Nagura\textquotesingle s theorem: for every real $x>25$, there is a prime in
$(x,6x/5)$; see~\cite{Nagura}. Choose $q<p<6q/5$, so $p\nmid q(q-1)$.
Since $d\le t<q$,
\[
 B(q,p)<\frac{1152}{25}q^{12}.
\]
The base point satisfies
\[
 25P(71)-1152\cdot71^{12}
 =79325064871285123771842993>0.
\]
If $q'$ is the next prime after $q$, Nagura\textquotesingle s theorem gives $q'<6q/5$,
and $P(q')=qP(q)$. Hence
\[
 \frac{P(q')}{(q')^{12}}>
 \frac{P(q)}{q^{12}}q(5/6)^{12}>
 \frac{P(q)}{q^{12}},
\]
where the last inequality follows from
$71\cdot5^{12}-6^{12}=15157202039>0$.
Induction along consecutive primes gives, for all $q\ge71$,
$P(q)>(1152/25)q^{12}>B(q,p)$.
Now $N\ge P(q)$ and Proposition~\ref{exp:B} complete the proof.
Finite verification is used only at the base points; the infinite range follows from the induction above.
\end{proof}

\subsection{Complete enumeration and verification of the remaining integers}

For $N\ge100000$, we have now covered the small-prime branches, all $q\ge61$, and
the cases $q\in\mathcal Q$ with $N\ge T(q)$.
The remainder is a known finite set: there are only finitely many $q$, and each has a finite interval.
Since $q=q(N)$ also implies $P(q)\mid N$, the candidates can be enumerated exactly as multiples.
Thus all remaining integers are precisely those of the form
\begin{equation}\label{exp:candidates}
 N=P(q)k,\quad 100000\le N<T(q),\quad q\in\mathcal Q,\quad q\nmid k.
\end{equation}
Conversely, these conditions ensure that every odd prime less than $q$ divides $N$, while $q\nmid N$,
so the least missing odd prime is indeed $q$. Distinct values of $q$ produce no duplicate candidates.

\begin{center}
\small
\begin{tabular}{rrrrr}
\toprule
$q$ & $d$ & $\ell$ & $T(q)$ & Candidates\\
\midrule
11&10&3&15488000&133229\\
13&6&5&18690048&14857\\
17&2&3&18939904&1181\\
19&6&5&89828352&333\\
23&22&3&720997376&142\\
29&14&3&1181544448&10\\
31&30&7&3321216000&1\\
37&6&5&1362604032&0\\
41&10&3&3442688000&0\\
43&42&5&17534555136&0\\
47&46&3&27521948672&0\\
53&26&3&25277943808&0\\
59&58&3&86935663616&0\\
\bottomrule
\end{tabular}
\end{center}

Put
\[
 L_q=\left\lceil\frac{100000}{P(q)}\right\rceil,\qquad
 U_q=\left\lfloor\frac{T(q)-1}{P(q)}\right\rfloor.
\]
If $L_q\le U_q$, the exact number of candidates is
\begin{equation}\label{exp:count}
 U_q-L_q+1-\left\lfloor U_q/q\right\rfloor
                  +\left\lfloor(L_q-1)/q\right\rfloor;
\end{equation}
and otherwise it is zero. Using $T(q)-1$ preserves the strict upper endpoint exactly.
The integer intervals for the nonempty rows are
\[
\begin{array}{c|r|r|r}
q&P(q)&L_q&U_q\\ \hline
11&105&953&147504\\
13&1155&87&16181\\
17&15015&7&1261\\
19&255255&1&351\\
23&4849845&1&148\\
29&111546435&1&10\\
31&3234846615&1&1
\end{array}.
\]
All other rows directly satisfy $P(q)>T(q)$.
The total number of candidates is
\[
 133229+14857+1181+333+142+10+1=149753.
\]
A valid common upper bound for the nonempty rows is $3321216000$,
and the largest actual candidate is $3234846615$.

For every remaining candidate, the witness table records three first summands $a_{--},a_{++},a_{-+}$.
Each, together with $N-a$, gives a concrete ordered decomposition, checked against the respective conditions
\[
 0<a<N,\qquad
 (g(a),g(N-a))=(-1,-1),(1,1),(-1,1).
\]
Exchanging the summands in the third decomposition gives $(1,-1)$.
All $449259$ concrete decompositions have been verified exactly.
The candidate set was also reconstructed independently starting from all multiples of $105$ and compared entry by entry with the certificate;
see Appendix~\ref{cert:complete}.

\subsection{Finite-range verification and completion of the proof}

Finally, consider $2\le N<100000$. Instead of scaling, we supply decomposition witnesses for each integer
and exhaust all decompositions for each genuinely absent pattern.
The finite witness table lists every integer in $2\le N\le99999$ consecutively, in $99998$ rows.
Each row records $(N,a_{--},a_{++},a_{-+})$. For every positive $a_{st}$, we independently check
\[
0<a_{st}<N,\qquad g(a_{st})=s,\qquad g(N-a_{st})=t.
\]
All $299981$ positive witnesses pass verification. The remaining $13$ positions are marked zero; for each, we exhaust
$1\le a<N$ to confirm absence of the specified pattern. Row-by-row checks also ensure consecutive indices, no duplicates, and no out-of-range entries.
The values $g(1),\ldots,g(99998)$ used in these checks are computed independently by complete integer trial division;
each distinct prime factor changes the sign once, and prime powers are not counted repeatedly.
Thus positive witnesses establish existence, and exhaustive checks establish absence; together they give the exact exceptional sets.
Exchanging the summands of a mixed-sign representation gives the other ordered mixed-sign pattern.
The three resulting zero sets are exactly those listed in Theorem~\ref{thm:main-explicit};
Appendix~\ref{app:small-range} gives the full verification conditions, data version, and all ordered counts for $2\le N\le17$.

We now combine the completed covering arguments and verifications.
\begin{proof}[Proof of Theorem~\ref{thm:main-explicit}]
First suppose $N\ge100000$, and take the least odd prime $q=q(N)$ not dividing $N$, as defined in~\eqref{glob:qdef}.
Then $N,q$ satisfy the relation $H(N,q)$.
\begin{enumerate}
\item If $q=3,5,7$, the small-prime branches of this section give all four patterns.
\item If $q\ge61$, Proposition~\ref{exp:tail} gives all four patterns.
\item Otherwise, $q\in\mathcal Q$. If $N\ge T(q)$, apply Proposition~\ref{exp:middle};
if $N<T(q)$, equation~\eqref{exp:candidates} places $N$ in the complete finite candidate set, and its three verified decompositions, together with reversal of the mixed pair, give all four patterns.
\end{enumerate}
These branches cover every $N\ge100000$; the threshold endpoint $N=T(q)$ belongs to the branch $N\ge T(q)$.

Now suppose $2\le N<100000$. The complete finite verification above determines the zero sets,
which are exactly the three exceptional sets in Theorem~\ref{thm:main-explicit}.
The ranges $2\le N<100000$ and $N\ge100000$ meet at consecutive integer endpoints and together exhaust all $N\ge2$. Hence these are the complete global exceptional sets.
Taking their union gives the necessary and sufficient condition for all four patterns to occur simultaneously.
The largest exceptions for the respective patterns are $3,7,17$, so their optimal eventual thresholds are $4,8,18$;
in particular, $R_{++}(17)=0$, whereas all four patterns occur for every $N\ge18$. Hence the common threshold $18$ is optimal.
\end{proof}

\phantomsection
\addcontentsline{toc}{section}{References}




\appendix


\section{Finite verification data and independent verification}\label{app:certificates}
The appendices describe the data formats and verification conditions and include selected witness tables; they do not reproduce every certificate entry. The complete certificates, witness data, and related verification programs are provided as electronic supplementary material in the GitHub repository:
\url{https://github.com/ilsunshine/Generalized_Goldbach_conjecture}.
The files described in this appendix are in the certificate directory of these materials.
For this version, the finite data were regenerated and checked entry by entry by an independent verifier that does not import the generation program.
The certificate version is determined by the fixed files and hash manifest in that directory.

\subsection{Data formats and completeness conditions}\label{cert:complete}
\begin{enumerate}
\item \texttt{chains.json} contains all covering chains. Each record has fields
\texttt{tag, q, d, ell, T, end, h, r, sign, j}.
The integer chain represented by a record is $n_i=d\ell j_i$. We check
$\ell j_i\equiv r\pmod h$, $g(d\ell j_i)=\varepsilon$, strict increase,
and $tn_0<T<qn_0$, $tn_{i+1}<qn_i$, $qn_J>\mathrm{end}$ individually.
For the intermediate-prime branches, $h=tq$, $T=128dq^2t^2$, and $\mathrm{end}=\ell T$;
the indices are exactly all $q\in\mathcal Q$, $0\le r<tq$ with $q\nmid r$, and both signs.
There are $31672$ intermediate chains containing $1770208$ coefficients.
\item The small-prime chains are in the same file. For $q=3,5$, use $\ell=5,3$, respectively,
covering $[100000,\ell\cdot100000]$.
The first group for $q=7$ uses $\ell=1$ and covers $[100000,1000000]$;
this group is used only for interval coverage, not for scaling.
The second group uses $\ell=5$, covers $[1000000,5000000]$, and then applies the covering-chain lemma.
Each group is indexed by $0\le r<tq$ with $q\nmid r$ and both signs,
corresponding to classes $dr$ modulo $tqd$. The small-prime chains total $184$ chains and $2203$ coefficients.
\item Each row of \texttt{odd.json} is $(q,d,r,\varepsilon,k)$, representing $A=dk$.
We check $k>0$, $k\equiv r\pmod{q^2}$, $g(dk)=\varepsilon$, and the complete class index.
For intermediate primes, we also check $k<128q^2$; there are $32492$ witnesses.
There are another $136$ witnesses for $q=3,5,7$; their strict upper bounds and the three supplementary decompositions for $q=7$ are checked as described in the main text.
\item Each row of \texttt{remaining.csv} is $(N,q,a_{--},a_{++},a_{-+})$.
Using~\eqref{exp:candidates}, we reconstruct the candidate set from multiples of each $P(q)$, check that it agrees exactly with the file\textquotesingle s
$(N,q)$ set without duplicates, and then check the strict intervals and signs for all summands.
There are $149753$ rows and $449259$ positive witnesses; the largest candidate is $3234846615$.
\end{enumerate}
Complete indices for covering chains, odd-branch witnesses, and remaining candidates are checked by set equality, not merely by comparing totals.
Completeness of this version follows from the necessary and sufficient conditions in~\eqref{exp:candidates}, explicit interval enumeration, and set equality.

\subsection{Sign computation and program independence}
\texttt{rebuild.py} generates the data; \texttt{verify.py} reads only the fixed data and does not import the generation program.
The generator uses a sieve that flips signs at prime multiples. The verifier independently constructs a smallest-prime-factor table up to $26000000$,
successively extracts each smallest prime factor, removes its full power, and computes the sign from the size of the prime-factor set.
For $cj$, it uses the union of prime-factor sets, avoiding any use of multiplicativity without a coprimality check.
Numbers beyond the table range are factored exactly by trial division by primes. The largest candidate is less than $3.33\cdot10^9$,
and the trial-division list contains every prime below $60000$, exceeding its square root.
Signs below one hundred thousand are computed separately by complete trial division for each integer, independently of the smallest-prime-factor table.
Endpoints, prime-prefix products, and base points use Python arbitrary-precision integers; all covering comparisons use integer arithmetic.

The smallest-prime-factor sieve is correct because every composite integer $n$ has a smallest prime factor $p\le\sqrt n$
and $n\ge p^2$. Thus the sieve visits $n$ when processing multiples of $p$, and taking minima leaves exactly its smallest prime factor.
Primes are never written by the sieves for smaller primes. During factorization, each prime factor is counted once, so the resulting sign is exactly $g$.
For numbers outside the table, once $p^2$ exceeds the remaining factor, that factor, if greater than $1$, must be prime.

\subsection{Reproduction and version records}
In an environment with Python 3 and NumPy installed, run
\texttt{python verify.py} to verify the fixed certificates; run
\texttt{python rebuild.py} to regenerate the witnesses.
\texttt{verification\_report.json} records the verified counts and exceptional sets for this run;
\texttt{sha256.json} records the SHA-256 hashes of the files.
The hashes identify byte versions only. The mathematical basis is the entry-by-entry verification of conditions and the infinite-range proofs in the main text.
The verification covers all $1772411$ chain coefficients, $32628$ odd-branch witnesses,
$99998$ rows below one hundred thousand, all remaining candidates, and the four exact integer base-point checks for the tail.

\Needspace{28\baselineskip}
\subsection{Short witness tables}\label{cert:smalltables}

The following are two complete short tables for the odd branches with small auxiliary primes.
In each row, $A_-,A_+$ are positive even integers with negative and positive signs, respectively,
and are congruent to the residue $r$ at the start of the row.

\begin{center}
\begin{tabular}[t]{rrr}
\toprule
$r\bmod9$&$A_-$&$A_+$\\
\midrule
1&64&10\\
2&2&20\\
4&4&22\\
5&32&14\\
7&16&34\\
8&8&26\\
\bottomrule
\end{tabular}
\qquad
\begin{tabular}[t]{rrr}
\toprule
$r\bmod25$&$A_-$&$A_+$\\
\midrule
1&126&26\\
2&2&52\\
3&78&28\\
4&4&54\\
6&156&6\\
7&32&82\\
8&8&58\\
9&84&34\\
11&186&36\\
12&312&12\\
13&138&38\\
14&64&14\\
16&16&116\\
17&42&92\\
18&168&18\\
19&294&44\\
21&246&46\\
22&222&22\\
23&198&48\\
24&174&24\\
\bottomrule
\end{tabular}
\end{center}

The longer integer chains and all remaining decompositions are included in the electronic supplementary material in the repository cited above.

\section{Complete finite verification below one hundred thousand}\label{app:small-range}

This appendix gives the complete finite verification for $2\le N<100000$, providing proof material for the finite-range part of Theorem~\ref{thm:main-explicit}. It checks the mathematical conditions of each row directly, without invoking the theorem\textquotesingle s conclusion to decide whether a witness is valid.

\subsection{Data and verification conditions}
The accompanying table \texttt{small\_range\_witnesses.csv} contains $99998$ rows, indexed consecutively by
\[
N=2,3,\ldots,99999.
\]
Each row contains four integers $(N,a_{--},a_{++},a_{-+})$. A positive $a_{st}$ is a witness; zero marks absence of that pattern. Each positive witness satisfies
\[
0<a_{st}<N,\qquad g(a_{st})=s,\qquad g(N-a_{st})=t.
\]
Reversing a negative-positive witness gives a positive-negative witness. All positive witnesses are at most $38$. This small bound has been checked only in the specified finite range; it does not assert that these thirty-eight anchors suffice for all integers.

The numbers of nonzero witnesses are, respectively,
\[
99996\quad(--),\qquad99989\quad(++),\qquad99996\quad(-+),
\]
for a total of $299981$. The independent verification performs complete trial division anew for all $1\le n\le99998$, obtains the parity of the number of distinct prime factors, and checks each witness individually. It does not import the sieve or sign-computation functions used to generate the witness table.

Trial division is justified as follows. Test successive divisors starting at $2$. When a prime factor is found, remove its entire power before changing the sign once. Once the square of the trial divisor exceeds the remainder, any remainder greater than $1$ must be prime, so the sign changes once more. Each distinct prime factor is counted exactly once, and the sign at $n=1$ is $1$.

For a zero mark, all $1\le a<N$ are checked separately to verify absence of the corresponding pattern; absence is not inferred from the anchor bound. All zero marks occur at $N\le17$ and can also be checked against the complete counts in the next subsection. The full row index, four-integer format of each row, strict intervals, and all sign conditions pass verification.

\subsection{Ordered counts for all small exceptions}
For $2\le N\le17$, the exact counts are as follows. Each row sums to $N-1$.
\begin{center}
\begin{tabular}{rrrrr}
\toprule
$N$&$R_{--}(N)$&$R_{-+}(N)$&$R_{+-}(N)$&$R_{++}(N)$\\
\midrule
2&0&0&0&1\\
3&0&1&1&0\\
4&1&1&1&0\\
5&2&1&1&0\\
6&3&1&1&0\\
7&4&0&0&2\\
8&3&2&2&0\\
9&4&2&2&0\\
10&5&2&2&0\\
11&6&1&1&2\\
12&6&2&2&1\\
13&6&2&2&2\\
14&5&4&4&0\\
15&6&3&3&2\\
16&7&2&2&4\\
17&4&6&6&0\\
\bottomrule
\end{tabular}
\end{center}
This table can also be checked directly by summing from $g(1),\ldots,g(16)$. For example, $2=1+1$ is a positive-positive representation, so $2$ is not an exception for that pattern.

For $18\le N\le99999$, all three witnesses in every row are positive, so all four patterns occur. Combining this with the argument for all $N\ge100000$ in Section~\ref{sec:explicit} gives the complete classification in Theorem~\ref{thm:main-explicit}.

\subsection{Data version and scope of the conclusions}
The decomposition table for this version is \texttt{small\_range\_witnesses.csv} in the accompanying certificate directory;
the hashes and actual verification records are in \texttt{sha256.json} and
\texttt{verification\_report.json}.
This appendix establishes existence and absence in the finite range.
The infinite range is handled by the covering chains and mathematical tail proof in Section~\ref{sec:explicit},
while the quantitative lower bounds follow from the separate argument in Section~\ref{sec:qplus}.

\section*{Statement on the use of tools}
The artificial intelligence tool GPT-6.1 sol was used to assist with organizing derivations, editing the text, and generating and verifying finite certificates.
\end{document}